\documentclass[a4paper,11pt,english]{article}
\usepackage{multirow}%
\usepackage{xcolor}%
\usepackage{fullpage}
\usepackage{amsmath}%
\usepackage{amsthm}%
\usepackage{amsfonts}%
\usepackage{amssymb}%
\usepackage{amsfonts}
\usepackage{graphicx}
\usepackage{colortbl}
\usepackage{enumerate}
\usepackage{dsfont}
\usepackage{titletoc}
\usepackage[runin]{abstract}
\usepackage{calc}

\newtheorem{theorem}{Theorem}

\newtheorem{proposition}[theorem]{Proposition}
\newtheorem{lemma}[theorem]{Lemma}

\newtheorem{remark}{Remark}
\newtheorem{definition}{Definition}
\newtheorem{example}{Example}

\newcommand{\indep}{\perp\hspace{-.25cm}\perp}
\newcommand{\E}{\mathbb{E}}
\renewcommand{\P}{\mathbb{P}}
\newcommand{\R}{\mathbb{R}}
\newcommand{\N}{\mathbb{N}}

\renewcommand{\^}[1]{\widehat{#1}}

\newcommand{\vecdeux}[2]{ \begin{pmatrix} #1 \\ #2 \end{pmatrix} }

\newcommand{\dd}{ \mathrm{d} }
\renewcommand{\r}{ \rightarrow }

\DeclareMathOperator{\im}{Im} 

\newcommand{\lr}{ \longrightarrow }

\DeclareMathOperator{\diag}{diag} 
 
\DeclareMathOperator{\var}{var} 
\DeclareMathOperator{\cov}{cov} 
\DeclareMathOperator{\rank}{rank} 
\DeclareMathOperator{\vecv}{vec} 
\DeclareMathOperator{\as}{a.s.} 
 
\DeclareMathOperator{\argmin}{argmin} 
\newcommand{\cRM}[1]{\MakeUppercase{\romannumeral #1}}
\newcommand{\cRm}[1]{\textsc{\romannumeral #1}}

\newcolumntype{C}[1]{>{\centering\arraybackslash }b{#1}}
\begin{document}

\title{Bootstrap Testing of the Rank of a Matrix via Least Squared
Constrained Estimation.}
\author{Fran\c cois Portier and Bernard Delyon}
\date{}
\maketitle

 \begin{abstract}
In order to test if an unknown matrix has a given rank (null hypothesis), we
consider the family of statistics that are minimum squared distances between an estimator
and the manifold of fixed-rank matrix. Under the null hypothesis, every statistic of this
family converges to a weighted chi-squared distribution. In this paper, we introduce the
constrained bootstrap to build bootstrap estimate of the law under the null hypothesis of
such statistics. As a result, the constrained bootstrap is employed to estimate the quantile
for testing the rank. We provide the consistency of the procedure and the simulations
shed light one the accuracy of the constrained bootstrap with respect to the traditional asymptotic comparison. More generally, the results are extended to test if an unknown parameter belongs to a sub-manifold locally smooth. Finally, the constrained bootstrap is easy to compute, it handles a large family of tests and it works under mild assumptions. 
 \bigskip
 
 \end{abstract}

\section{Introduction}\label{s1}

Let $M_0\in \R^{p\times H} $ be an unknown matrix (arbitrarily $p \leq H$). To infer about the rank of $M_0$ with hypothesis testing, the general framework usually considered is the following: there exists an estimator $\^M\in \R^{p\times H}$ of $M_0$ such that
\begin{align}\label{as1}
n^{1/2}(\^M-M_0) \overset{\dd}{\lr} W, \qquad \text{with} \qquad \vecv{(W)} = \mathcal N(0,\Gamma) 
\end{align}
where $\vecv(\cdot)$ vectorizes a matrix by stacking its columns. In the whole paper the hatted quantities are random sequences that depends on the sample number $n$, all the limit are taken with respect to $n$. Moreover there exists an estimator $\^\Gamma$ such that
\begin{align}\label{as2}
\^\Gamma \overset{\P}{\lr} \Gamma, 
\end{align}
and in some cases, one may ask that
\begin{align}\label{as3}
\Gamma \text{ is full rank.} 
\end{align}
Let $d_0$ be the rank of $M_0$ and $m \in \{1, ..., p\}$, we consider the set of hypotheses
\begin{align}\label{hyp}
H_0\ :\quad d_0 = m \qquad \text{against}\qquad H_1\ :\quad d_0 > m,
\end{align}
Thus $d_0$ can be estimated the following way: we start by testing $m = 0$, if $H_0$ is rejected we go a step further $m := m + 1$, if not we stop the procedure and the estimated rank is $\^d = m$. In this paper, by considering the hypotheses (\ref{hyp}) we focus on each step of this procedure. 

Many different statistical tests appeared in the literature for this purpose. For instance Cragg and Donald \cite{donald1996} introduced a statistic based on the LU decomposition of $\^M$, Kleibergen and Paap \cite{paap2006} studied the asymptotic behaviour of some transformation of the singular values of $\^M$, and Cragg and Donald \cite{donald1997} considered the minimum of a
squared distance under rank constraint. In some other fields with similar issues,
close ideas have been developed : Bura and Yang \cite{bura2010} examined a Wald type statistic depending on the singular decomposition of $\^M$ and Cook and Ni \cite{cook2005} also considered the minimum of a squared distance under rank constraint. Although based
on different considerations, each of the previous work relies on the test described by (\ref{hyp}). For
comprehensiveness, in this paper we consider the following three statistics. The first one is
introduced by Li \cite{li1991} as
\begin{align}\label{stat1}
\^\Lambda_1 = n\sum_{k=m+1}^p \^\lambda_k^2
\end{align}
where $(\^\lambda_1, ..., \^\lambda_p)$ are the singular values of $\^M$ arranged in descending order. Under $H_0$ and
(\ref{as1}), this statistic converges in law to a weighted chi-squared distribution \cite{bura2010}. The main drawback of such a test is that $\^\Lambda_1$ is not pivotal, i.e. its asymptotic law depends on
unknown quantities that are $M_0$ and $\Gamma$. Accordingly the consistency of the associated test requires assumptions
(\ref{as1}) and (\ref{as2}). In \cite{bura2010} a standardized version of $\^\Lambda_1$ is studied with
\begin{align}\label{stat2}
\^\Lambda_1 = n \vecv(\^Q_1\^M\^Q_2)^T [(\^Q_2 \otimes \^Q_1)\^\Gamma (\^Q_2 \otimes \^Q_1)]^+ \vecv(\^Q1\^M\^Q2)
\end{align}
where $M^+$ stands for the Moore-Penrose inverse of $M$ and $\^Q_1$ and $\^Q_2$ are respectively the orthogonal projectors on the left and right singular spaces associated with the $p-m$ smallest singular
values of $\^M$. The authors proved that under $H_0$, if (\ref{as1}) and (\ref{as2}) hold, the Wald-type statistic $\^\Lambda_2$ is asymptotically chi-squared distributed. Besides, \cite{donald1997} and \cite{cook2005} proposed a constrained estimator
by minimizing a squared distance under a fixed-rank constraint as
\begin{align}\label{stat3}
\^\Lambda_3 = n\underset{\rank(M) = m}{\min}  \vecv(\^M-M)^T \^\Gamma^{-1} \vecv(\^M-M),
\end{align}
which is also asymptotically chi-squared distributed under $H_0$, 
assuming (\ref{as1}), (\ref{as2}) and (\ref{as3}). We
will refer the minimum discrepancy approach. Although the statistics $\^\Lambda_2$ and $\^\Lambda_3$ have the
convenience of being pivotal, they both require the inversion of a large matrix and this may
cause robustness problems when the sample number is not large enough. For $\alpha\in ]0,1[$ and
under the relevant assumptions, each of these statistics $\^\Lambda_1$, $\^\Lambda_2$ and $\^\Lambda_3$, is consistent at level $\alpha$ in testing (\ref{hyp}), i.e. the level goes to $1-\alpha$ and the power goes to $1$ as $n$ goes to $\infty$. 

Nevertheless the estimation of the quantile is difficult because either the asymptotic distribution depends on the data (non pivotality represented by $\^\Lambda_1$), or the true distribution may be quite different than the asymptotic one (slow rates of convergence represented by $\^\Lambda_2$ and $\^\Lambda_3$). The objective of the
paper is to propose a bootstrap method for quantile estimation in this context. 

An important remark which instigates the sketch of the paper is that all the previous statistics share the form
\begin{align}\label{sharedform}
\^\Lambda = n  \|\^B \vecv(\^M-\^M_c)\|^2  \quad \text{with} \quad \^M_c = \underset{\rank(M)=m}{\argmin}  \|\^A\vecv(\^M-M)\|^2 
\end{align}
where $\|\cdot\| $ is the Euclidean norm, $\^A\in \R^{pH\times pH}$, $\^B\in \R^{pH\times pH}$. The values of $\^A$ and $\^B$ corresponding to the statistics $\^\Lambda_1$, $\^\Lambda_2$ and $\^\Lambda_3$ are summarized in the Table \ref{table1} (See Section
\ref{s2} for the details).

\begin{table}
  \centering
\begin{tabular}{ c | C{4cm} C{4cm} C{4cm}}
   & \multicolumn{1}{ >{\columncolor[gray]{.7}}c }{$\^\Lambda_1$}& \multicolumn{1}{ >{\columncolor[gray]{.7}}c }{$\^\Lambda_2$} & \multicolumn{1}{ >{\columncolor[gray]{.7}}c }{$\^\Lambda_3$} \\ 
 \hline
  \hline
\multirow{2}{*}{$\^A$} & \multirow{2}{*}{$I$}& \multirow{2}{*}{$I$}& \multirow{2}{*}{$\^\Gamma^{-1/2}$} \\ 
 &  &  & \\  \hline
\multirow{2}{*}{$\^B$}& \multirow{2}{*}{$I$}& \multirow{2}{*}{$[(\^Q_2\otimes \^Q_1)\^\Gamma (\^Q_2\otimes \^Q_1)]^{+1/2}$}& \multirow{2}{*}{$\^\Gamma^{-1/2}$} \\ 
 &  &  & \\  \hline
\end{tabular}
 \caption{Values of $\^A$ and $\^B$ in (\ref{sharedform}) for $\^\Lambda_1$, $\^\Lambda_2$ and $\^\Lambda_3$.}\label{table1}
\end{table}

We refer to traditional testing (resp. bootstrap testing) when the statistic is compared to
its asymptotic quantile (resp. bootstrap quantile). The bootstrap test is said to be consistent
at level $\alpha$ if
\begin{align}\label{levelp}
\P_{H_0}\left(\^\Lambda >\^q(\alpha)\right) \overset{}{\lr}1-\alpha
\qquad \text{and} \qquad \P_{H_1}\left(\^\Lambda 
>\^q(\alpha)\right)\overset{}{\lr} 1,
\end{align}
where $\^q(\alpha)$ is the quantile of level $\alpha$ calculated by bootstrap. The advantage of bootstrap testing is its high level of accuracy under $H_0$ with respect to traditional testing. This fact is emphasized by considering the two possibilities: when the statistic is pivotal and when the asymptotic law of the statistic depends on unknown quantities. First, as highlighted by Hall \cite{hall1992}, when the statistic is pivotal, under some conditions the gap between the distribution of the statistic and its bootstrap distribution is $O_{\P}(n^{-1})$. Since the normal approximation leads to a difference $O(n^{−1/2})$, the bootstrap enjoys a better level of accuracy. Secondly if the asymptotic law of the statistic is unknown, the bootstrap appears even more as a convenient alternative because it avoids its estimation. In \cite{hall1991}, Hall and Wilson give two advices for the use of the bootstrap testing:
\begin{enumerate}[A)]
\item  \label{gui1}Whatever the sample is under $H_0$ or $H_1$, the bootstrap estimates the law of the statistic
under $H_0$.
\item \label{gui2}The statistic is pivotal.
\end{enumerate}
The first guideline is the most crucial because if it fails it may lead to inconsistency of the test.
The second guideline aims at improving the accuracy of the test by taking full advantage of the
accuracy of the bootstrap. In this paper we propose a new procedure for bootstrap testing in least square constraint estimation (LSCE) (estimators as (\ref{sharedform}) are particular cases), called constrained bootstrap (CS bootstrap). More precisely, the CS bootstrap aims at testing whether a parameter belongs or not to a submanifold and so generalised the test (\ref{hyp}). Our main result is the consistency of the CS bootstrap under mild conditions. As a consequence we provide a consistent bootstrap testing procedure for testing (\ref{hyp}) with the statistic $\^\Lambda_1$, $\^\Lambda_2$ and $\^\Lambda_3$. For the sake of clarity, we address
the CS bootstrap in the next section. Section 3 is dedicated to rank estimation with special
interest to the bootstrap of the statistic $\^\Lambda_1$, $\^\Lambda_2$ and $\^\Lambda_3$. Finally, the last section emphasizes the accuracy of the bootstrap in rank estimation by providing a simulation study in sufficient
dimension reduction (SDR). Accordingly, the sketch of the paper is as follows:

$\bullet$ The CS bootstrap in LSCE

$\bullet$ Bootstrap testing procedure for $\^\Lambda_1$, $\^\Lambda_2$ and $\^\Lambda_3$

$\bullet$ Application to SDR

\section{The constrained bootstrap for LSCE and hypothesis testing}\label{s2}

Because of (\ref{sharedform}) LSCE has a central place in the paper. Moreover since LSCE intervenes in many statistical fields as M-estimation or hypothesis testing, this section is independent from the rest of the paper.

\subsection{LSCE}\label{s21}
Let $\theta_0 \in\R^p$ be called the parameter of interest, and let $\^\theta\in \R^p$ be an estimator of $\theta_0$. We define the
constrained estimator of $\theta_0$ as
\begin{align}\label{estimatorconst}
\^\theta_c = \underset{\theta\in \mathcal M}{\argmin}\ (\^\theta- \theta)^T \^A (\^\theta - \theta) , 
\end{align}
where $\mathcal M$ is a submanifold of $\R^p$ with co-dimension $q$, and  $\^A \in \R^{p\times p}$. The constrained statistic is defined as
\begin{align}\label{minchi2}
\^\Lambda =n (\^\theta- \^\theta_c)^T \^B (\^\theta - 
\^\theta_c).
\end{align}
where $\^B\in \R^{p\times p}$. Note that if $\^A$ is full rank, the unique minimizer of (\ref{estimatorconst}) without constraint is $\^\theta$, hence it could be understood as the unconstrained estimator. We introduce now the notion of nonsingular point in $\mathcal M$. This one is needed to express the Lagrangian first order condition of the optimization (\ref{estimatorconst}). For any 
function $g=(g_1,\ldots,g_p):\R^p\r \R^q$, define its Jacobian as $J_g 
=(\nabla g_1,..., \nabla g_{q})$, where $\nabla$ stands for the gradient operator.

\begin{definition}\label{const}
We say that $\theta$ is $\mathcal M$-nonsingular if $\theta\in \mathcal M$ and if there 
exists a neighbourhood $V$ and a function $g:\R^p\r \R^q$ continuously differentiable on $V$ with $J_g(\theta)$ full rank such that
\begin{align*}
V\cap \mathcal M =\{g=0\}.
\end{align*}
\end{definition}
As a consequence any point of a submanifold locally smooth is nonsingular, e.g. any matrix with rank $m$ is a nonsingular point in the submanifold $\rank(M)= m$. We prove in Proposition \ref{ssboot2} that if $\theta_0$ is $\mathcal 
M$-nonsingular, $\sqrt n (\^\theta-\theta_0)\overset{\dd}{\r} 
\mathcal N (0,\Delta)$ and $\^B=\^A\overset{\P}{\r} A$ is full rank, then we have
\begin{align}\label{ssboot}
\^\Lambda\overset{ \dd } {\lr } \sum_{k=1}^p \nu_k W_k^2,   
\end{align}
where the $W_k$'s are i.i.d. Gaussian random variables and the $\nu_k$'s are the eigenvalues of the matrix $\Delta^{1/2} J_g(\theta_0)^T(J_g(\theta_0)A^{-1}J_g(\theta_0)^T)^{-1} J_g(\theta_0)\Delta^{1/2}$.  Especially, the case $A=\Delta^{-1}$ is interesting because $\^\Lambda$ is asymptotically chi-squared distributed with $q$ degrees of freedom. 
Otherwise, if $\theta_0\notin \mathcal M$, $\^\Lambda$ goes to infinity in probability. Those facts shed light on a consistent testing procedure based on LSCE with the hypotheses
\begin{align}\label{h0}
H_0\ : \qquad \theta_0\in \mathcal M  \qquad \text{against} \qquad 
H_1\ : \quad \theta_0\notin \mathcal M 
\end{align}
and the decision rule to reject $H_0$ if $\^\Lambda$ is larger than 
a quantile of its asymptotic law. Accordingly the previous framework can be seen as an extension of the Wald test statistic which handles the simple hypothesis $\theta_0=\theta$ with the statistic $(\^\theta - \theta )^T \Delta^{-1} (\^\theta - \theta )$.

\subsection{The bootstrap in LSCE}

Since LSCE is a particular case of estimating equation, we review the bootstrap literature with two principal directions: estimating equation and hypothesis testing. For clarity we alleviate the framework in this section: let $X_1,\cdots, X_n$ be an i.i.d. sequence of real random variables with law $P$, define $\gamma = \var(X_1)$, $\^\gamma = \overline{(X - \overline{X})^2} $, we put $\theta_0 = \E[X_1]$, $\^\theta = \overline{X}$, and  $A = B = \gamma^{-1 }$ where $\overline{\ \cdot\ }$ stands
for the empirical mean.

The original bootstrap was introduced in \cite{efron1979} in the following way. Let $X_1^*,\ldots, X_n^*$ be an i.i.d. sequence of real random variables with law $\^P = n^{-1}\sum_{i=1}^n \delta_{X_i}$, define $\theta^*=\overline{X^*}$, the distribution of $\sqrt n (\theta^*-\^\theta)$ conditionally on the sample, that we call the bootstrap distribution, is ``close" to the distribution of $\sqrt n ( \^\theta-\theta_0)$, that we call the true distribution (in the rest of the paper we just say ``conditionally" instead of ``conditionally on the sample"). For instance, it is shown in \cite{vandervaart1998} that the bootstrap distribution converges weakly to the true distribution almost surely. One says that $\sqrt n (\theta^*-\^\theta)$ bootstraps $\sqrt n ( \^\theta-\theta_0)$ and we will write 
\begin{align*}
\mathcal L_{\infty} (n^{1/2} (\theta^*-\^\theta)|\^P)= \mathcal L_{\infty} (n^{1/2} ( \^\theta-\theta_0) )\quad \as,
\end{align*} 
where $\mathcal L_{\infty}(\cdot)$ and $\mathcal L_{\infty}(\cdot|\^P)$ both mean the asymptotic laws with the difference that the later is conditional on the sample. Equivalently, one has for every $x\in \R$, $\P(\sqrt n  (\theta^*-\^\theta)\leq x |\^P) \overset{\as}{\r} \P(\sqrt n ( \^\theta-\theta_0)\leq x )$, but the use of the bootstrap is legitimate by a more general results stated in \cite{hall1992}, which says that
\begin{align}\label{hall1}
|\P(n^{1/2} (\theta^*-\^\theta)/\gamma^*\leq x |\^P) - \P(n^{1/2} ( \^\theta-\theta_0)/\^\gamma\leq x )| = O_{\P}(n^{-1})
\end{align}
with $\gamma^* = \overline{(X^*-\overline{X^*})^2 }$, provided that $P$ is non-lattice. Besides, one has
\begin{align*}
|\P(n^{1/2}( \^\theta-\theta_0)/\^\gamma\leq x )-\Phi(x) |= O_{\P}(n^{-1/2}),
\end{align*}
where $\Phi$ is the cumulative distribution function (c.d.f.) of the standard normal law. Variations of Efron's resampling plan are proposed in \cite{barbe1995} under the name of weighted bootstrap. For a complete introduction about the bootstrap we refer to \cite{hall1992}. We now present three different bootstrap techniques related to LSCE\footnote{A bootstrap with a Delta-method approach (see \cite{vandervaart1998}, chapter 23, Theorem 5) fails because $x\r \underset{\|\theta\|=1}{\min} \|x  - \theta \|$ is not 
continuously differentiable on the unit circle.}.

\bigskip

\begin{enumerate}[(i)]
\item \label{cboot} The classical bootstrap (C bootstrap)

The literature about the bootstrap in Z and M-estimation, see respectively \cite{chatterjee2005} and \cite{arcones1990}, is based on the following principle: if $\theta_M=\underset{\theta\in \Theta}{\argmin}\ \E[ \phi(X,\theta)]$ is estimated by $\^\theta_M = \underset{\theta\in \Theta}{\argmin}\ \frac 1 n 
\sum_{i=1}^n \phi(X_i,\theta)$ where $\Theta$ is an open set, then the bootstrap of $\sqrt n (\^\theta_M- \theta_M)$ is carried out by the quantity $\sqrt n (\theta_M^*- \^\theta_M)$ with 
\begin{align}\label{boot1}
\theta_M^* = \underset{\theta\in \Theta}{\text{argmin}} \ n^{-1} \sum_{i=1}^n w_i \phi(X_i,\theta) ,
\end{align}
where $(w_i)$ is a sequence of random variables. The particular case where the vector
$(w_1 ,\ldots ,w_n)$ is distributed as $\text{mult} (n,(n^{-1},\ldots,n^{-1}))$ leads to a direct 
application of original Efron's bootstrap to M-estimation. Since 
such a bootstrap has been extensively studied, we refer to the C bootstrap. To the knowledge of the authors, the C bootstrap when $\Theta$ has empty interior has not been studied 
yet. Nevertheless one may sight its bad behaviour for the test of equal mean $H_0:\ \theta_0 = \mu$. The associated least squared constrained statistic
\begin{align*}
n\^\gamma^{-1} (\^\theta - \mu)^2,
\end{align*}
is indeed the score statistic associated to the M-estimator with $\phi(x,\theta)=\^\gamma ^{-1}(x -\theta)^2$ and $\Theta = \{\mu\}$. Clearly the C bootstrap through $n \gamma^{*-1}  (\theta^*- \mu)^2$ does not work because of its bad behaviour under $H_1$ for instance. In this case it is better to use
\begin{align*}
n \gamma^{*-1} (\theta^* - \^\theta)^2,
\end{align*}
but it cannot handle the cases of more involved hypotheses\footnote{We refer to \cite{hall1991} for a study of this bootstrap in order to test $\theta_0=\mu$.}. Whereas the C bootstrap is not really connected with hypothesis testing, the two following bootstrap procedures are more related to the present work.

\item \label{bboot} The biaised bootstrap (B bootstrap)

The B bootstrap is introduced in \cite{hall1999} and is directly motivated by hypothesis testing. The original idea of their work is to re-sample with respect to the distribution $\^P_b = \sum_{i=1}^n \omega_i \delta_{X_i} $, where the $\omega_i$'s maximize 
\begin{align}\label{hallboot}
\sum_{i=1}^n \log(\omega_i) \quad \text{under the constraints}\quad \begin{array}{r} \frac 1 n \sum_{i=1}^n \omega_i X_i= \mu \\ \sum_{i=1}^n \omega_i =1 \end{array}.
\end{align}
Since the $\omega_i$'s minimize the Kulback-Leibler distance 
between $\^P$ and $ \^P_b$, one can see the resulting distribution 
as the closest to the original one satisfying the mean constraint. 
The authors presented interesting results for the test of equal 
mean $\theta_0 = \mu$, essentially the bootstrap statistic $ n 
\gamma^{*-1}( \theta_b^*-\mu)^2$, with $\theta_b^*=\overline{X^*_b}$, 
$X_{b,i}^*$ sampled from $\^P_b$, has a chi-squared limiting distribution either 
$H_0$ or $H_1$ is assumed. As a result both guidelines (\ref{gui1}) 
and (\ref{gui2}) are checked. They go further by showing that the B 
bootstrap outclasses the asymptotic normal approximation for quantile 
estimation in the sense that 
$|\^q(\alpha)-q_{n}(\alpha)|=O_{\P}(n^{-1})$ whereas $|q_{n}(\alpha)- 
q_{\infty}(\alpha)|=O_{}(n^{-1/2})$, where $q_{\infty}$, $q_{n}$ and 
$\^q_n$ are the quantile functions of the standard normal distribution, 
the statistic $ n \^\gamma^{-1}(\^\theta - \mu )^2$ under $H_0$ 
and the bootstrapped statistic, respectively. Although the B bootstrap matches the context 
of hypothesis testing, it has been designed to handle the particular 
test of equal mean. To the knowledge of the authors the study of the B 
bootstrap has not been extended to other tests. Facing 
(\ref{hallboot}), the main drawback of the B bootstrap deals with 
algorithmic difficulties. Indeed when the constraint becomes more 
involved, solving (\ref{hallboot}) is more 
difficult. As a result it is not sure that this method could 
handle other situations such as fixed-rank constraints.

\item \label{efboot} The estimating function bootstrap (EF bootstrap)

Now $X_i\in \R^p$. Some other ideas about the 
bootstrap of the $Z$-estimators can be found in \cite{lele1991} and 
\cite{hu2000}, and can be summarized as follows. Considering the score statistic $\^S=\sqrt n \sum_{i=1}^n \frac{\partial \phi}{\partial \theta} (X_i,\theta_0)$, \cite{hu2000} showed that it could be bootstrapped by 
\begin{align*}
 S^*=n^{-1/2} \sum_{i=1}^n w_i \frac{\partial \phi}{\partial \theta} (X_i, \^\theta),
\end{align*}
where $(w_i)$ is a sequence of random variables. This bootstrap is called the EF bootstrap
and revealed nice computational properties. Moreover the 
authors argued for its use in quantile estimation in order to test if 
$g(\theta_0)=0$, where $g:\R^p\r\R^q$ is the constraint function, by recommending essentially to use $S^{*T}  J_g(\^\theta)^T\left(J_g(\^\theta) \gamma^{*} J_g(\^\theta)^T\right)^{-1} J_g (\^\theta) S^*$. Applying it to the least squared context $\phi(x,\theta)=\| \^\gamma ^{-1/2}(x -\theta)\|^2$, the EF bootstrap is carried out by
\begin{align*}
n(\theta^* - \^\theta)^T J_g(\^\theta)^T\left(J_g(\^\theta) \gamma^* J_g(\^\theta)^T\right)^{-1} J_g (\^\theta)  (\theta^* - \^\theta).
\end{align*}
Although it verifies both guidelines (\ref{gui1}) and (\ref{gui2}) (see the article for details), one can see that the good behaviour of 
such an approach is more based on the rank deficiency of 
$J_g(\^\theta)$ than on the bootstrap of $\sqrt n (\^\theta- \^\theta_c)$. Indeed $\sqrt n (\theta^* - \^\theta)$ bootstraps the non constrained estimator $\sqrt n (\^\theta - \theta_0)$. Then as the authors noticed, it is 
first of all a bootstrap of the Wald-type statistic $n \^S^{T}  J_g(\theta_0)^T\left(J_g(\theta_0) \^\gamma J_g(\theta_0)^T\right)^{-1} J_g (\theta_0) \^S$ which has fortunately the same asymptotic law than the targeted one. This may induce some loss in accuracy. Moreover, it requires the knowledge of the function $J_g$ which is not the case for fixed rank constraints where the $g$ depends on the limit $M_0$ (see Remark \ref{ex1} for some details).
\end{enumerate}

Essentially both (\ref{cboot}) and (\ref{bboot}) provide a bootstrap for testing simple hypotheses. The EF bootstrap proposed in (\ref{efboot}) extends this limited scope by including tests of the form $g(\theta_0)=0$ where $g$ is known. Nevertheless it does not handle the test (\ref{hyp}) as it is highlighted by the following remark.

\begin{remark}\label{ex1} \normalfont
Testing (\ref{hyp}) with $\^\Lambda_3$ results in an optimization with the 
constraint $\rank(M)=m$. Since the subspace of fixed rank matrices is 
a submanifold locally smooth
 with co-dimension $(p-d)(H-d)$, at every point $M$, there 
exists a neighbourhood $V$ and a $\mathcal C ^{\infty}$ function $g:V\r 
\R^{(p-d)(H-d)} $ such that $V\cap 
\{\rank(M)= m\} =\{g= 0\}$ and $J_{g}(M)$ has full rank. Moreover, we have
\begin{align*}
\|\Gamma^{-1/2} \vecv(\^M_c-M_0)\| \leq  2 \|\Gamma^{-1/2} \vecv(\^M-M_0)\|.
\end{align*}
If now (\ref{as1}) holds, the right-hand side term goes to $0$ in 
probability and $\^M_c\overset{\P}{\r}M_0$. As a consequence, if $\Gamma$ is invertible, for any neighbourhood of $M_0$, from a certain rank, $\^M_c$ belongs to it with probability $1$. Then under $H_0$ since $M_0$ has rank $m$ the constrained estimator has the expression
\begin{align*}
\^M_c = \underset{g(M)= 0 }{\argmin} \|\Gamma^{-1/2} \vecv(\^M_c-M)\|,
\end{align*}
with $g$ depending on $M_0$. Unfortunately we do not know neither $g$ nor $J_{g}(M_0)$. This entails some problems relating to the later approach.
\end{remark}

\subsection{The constrained bootstrap}

The CS bootstrap is introduced in order to solve all the issues we have raised through the previous little review which are essentially: computational difficulties and small scope of the existing methods. The CS bootstrap targets an estimation $\^q(\alpha)$ of the quantile under $H_0$ of $\^\Lambda$. The consistency of the procedure, i.e. (\ref{levelp}), forms the main result about the CS bootstrap. Another important issue which occurs beforehand in the section is the bootstrap of the law of
\begin{align*}
n^{1/2} (\^\theta_c - \theta_0)\quad \text{under $H_0$}.
\end{align*}
Basically, we show that a bootstrap of the unconstrained estimator $\sqrt n (\^\theta- \theta_0)$ allows a bootstrap of the constrained estimator $\sqrt n (\^\theta_c - \theta_0)$ under $H_0$. We point out that the CS bootstrap heuristic is rather different 
than the C and EF bootstrap. Otherwise it shares the idea to ``reproduce" $H_0$ even if $H_1$ 
is realized with the B bootstrap. Assuming that we can bootstrap $\sqrt n (\^\theta- \theta_0)$, the CS bootstrap calculation of the statistic is realized as follows:

\noindent\fbox{\parbox{\linewidth-2\fboxrule-2\fboxsep}{
     \noindent\textbf{The CS bootstrap procedure} 

\noindent Compute 
\begin{align}\label{bootcond1}
\theta^*_0 = \^\theta_c + n^{-1/2}W^*, \qquad \text{with} \quad \mathcal L_{\infty}( W^*|\^P) = \mathcal L_{\infty} (n^{1/2} (\^\theta-\theta_0))\quad \as,
\end{align}
where the simulation of $W^*$ can be done by a standard bootstrap procedure\footnotemark. Calculate
\begin{align}\label{statboot}
\theta_c^* = \underset{\theta\in\mathcal M}{\argmin }\ (\theta^*_0- \theta)^T  A^{*} (\theta^*_0-\theta),\qquad \text{and} \qquad 
 \Lambda^* = n(\theta^*_0- \theta^*_c)^T  B^{*} (\theta^*_0-\theta^*_c),
\end{align}
where $A^*\in \R^{p\times p}$ and $B^*\in \R^{p\times p}$ \footnotemark.}}
\footnotetext[3]{The bootstrap procedure to get $W^*$ is not 
specified because it depends on $\^\theta$. For instance, if $\^\theta$ 
is a mean over some i.i.d. random variables, one can use the Efron's 
traditional bootstrap and if $\^\theta$ is a M-estimator, one should 
use a bootstrap as detailed by equation (\ref{boot1}).}
\footnotetext[4]{Assumptions about $A^*$ and $B^*$ are provided further in the statements of the propositions.}

Intuitively, this choice appears natural because $\theta^*_0$ equals 
$\^\theta_c$ plus a small perturbation going to $0$. Accordingly 
$\theta^*_0$ is somewhat reproducing the behaviour of $\^\theta$ under 
$H_0$, especially because $ W^*$ has the right asymptotic variance. As we should notice, $A^*$ and $B^*$ could be chosen as $\^A$ and $\^B$ but this is not the best choice in practice. As it is highlighted in (\ref{hall1}), we should normalized by the associated  bootstrap quantities (e.g. the variance computed on the bootstrap sample). The following lemma gives a first order 
decomposition of the bootstrap law $\sqrt n (\theta^*_c-\^\theta_c)$ 
under mild conditions. The following lemma is proved 
in the Appendix.

\begin{lemma}\label{bootth}
Let $\mathcal M$ be a submanifold. Assume there exists $\^\theta_c\in 
\mathcal M$ and $\theta_c$ a $\mathcal M$-nonsingular point such that
$\^\theta_c\overset{\as}{\r} \theta_c$. If moreover $\mathcal L _{\infty} (\sqrt 
n (\theta^*_0-\^\theta_c)|\^P)$ exists a.s. and conditionally a.s. $A^* 
\overset{\P}{\r}A$ is full rank, then we have conditionally a.s.
\begin{align*}
n^{1/2} (\theta_c^* -\^\theta_c) = (I-P) n^{1/2} (\theta^*_0-\^\theta_c)  +o_{\P}(1),
\end{align*}
with $P= A^{-1} J^T_g(\theta_c)(J_g(\theta_c) A^{-1} J^T_g(\theta_c))^{-1} J_g(\theta_c)$.
\end{lemma}

Note that if $\theta_0$ is $\mathcal M$-nonsingular and $\mathcal L _{\infty} (\sqrt 
n (\^\theta-\theta_0)|\^P)$ exists, we can apply Lemma \ref{bootth} with 
$\^\theta_c =  \theta_c= \theta_0$. This gives the following proposition:

\begin{proposition}\label{ssboot2}
Let $\mathcal M$ be a submanifold. Assume that $\mathcal L _{\infty} (\sqrt 
n (\^\theta-\theta_0)|\^P)$ exists with $\theta_0$  $\mathcal 
M$-nonsingular. Assume also that $\^A \overset{\P}{\r}A$ is full 
rank, then we have
\begin{align*}
n^{1/2} (\^\theta_c -\theta_0) = (I-P) n^{1/2} (\^\theta-\theta_0)  +o_{\P}(1),
\end{align*}
with $P= A^{-1} J^T_g(\theta_0)(J_g(\theta_0) A^{-1} 
J^T_g(\theta_0))^{-1} J_g(\theta_0)$.
\end{proposition}

Proposition \ref{ssboot2} leads easily to (\ref{ssboot}) and extends 
classical results  \cite{boos1990} about constrained estimators with constraint 
$\{g=0\}$ to manifold type constraints. Besides statements 
of Lemma \ref{bootth} and Proposition \ref{ssboot2} together explain 
the preceding definition of $\theta^*_0$ in (\ref{bootcond1}). They also 
lead to the following theorem.

\begin{theorem}\label{bootcor}
Let $\mathcal M$ be a submanifold. Assume that $\^\theta \overset{\as}{\r} \theta_0$ with $\theta_0$ $\mathcal M$-nonsingular and $\^A\overset{\P}{\r}A$ hold. If moreover (\ref{bootcond1}) holds and conditionally a.s. $A^* \overset{\P}{\r}A$ is full rank, then we have
\begin{align*}
\mathcal L_{\infty} (n^{1/2}(\theta^*_c-\^\theta_c)|\^P) = \mathcal 
L_{\infty} (n^{1/2} ( \^\theta_c - \theta_0  ))\quad \as .
\end{align*}

\end{theorem}

Essentially, Theorem \ref{bootcor} is an application of Lemma 
\ref{bootth} under $H_0$, indeed as we seen in the proof of Lemma \ref{bootth}, equation (\ref{ineg}), the assumption 
$\^\theta\overset{\as}{\r} \theta_0 \in \mathcal M$ implies that 
$\^\theta_c\overset{\as}{\r} \theta_c$. Nevertheless under $H_1$ nothing guarantee such a convergence (see Example \ref{ex2} below). Roughly speaking, asking for an equality in law under $H_1$ as in Theorem \ref{bootcor} may be too much to ask. However as stated in the following theorem we do not require that $\^\theta_c$ converges a.s. to a constant to provide that the power of the corresponding test goes to $1$. This leads to the consistency of the CS bootstrap for hypothesis testing. For the statement of the consistency theorem, we need to define the quantile function of the bootstrap statistic
\begin{align*}
\^q(\alpha) =  \text{inf} \ \{ x \ :\ \^F(x)\geq 1-\alpha\},
\end{align*}
where $\^F$ is the c.d.f. of $\Lambda^*$ 
conditionally on the sample. 

\begin{theorem}\label{bootthtest}
Let $\mathcal 
M$ be a manifold. Assume that $\^\theta\overset{\as}{\r} \theta_0$ 
with $\theta_0$ $\mathcal M$-nonsingular under $H_0$. We assume also 
that $\^A \overset{\P}{\r}A$ is full rank, $\^B\overset{\P}{\r}B$. If moreover $\mathcal L_{\infty} (\sqrt n (\theta^*_0-\^\theta_c)|\^P){=}\mathcal L_{\infty} (\sqrt n (\^\theta-\theta_0))  $ a.s. has a density, and conditionally a.s. $A^* \overset{\P}{\r}A$, $B^* \overset{\P}{\r} B$, then we have
\begin{align*}
\P_{H_0}(\^\Lambda >  \^q(\alpha) )\overset{}{\lr}1-\alpha, \quad \text{and} \quad 
\P_{H_1}(\^\Lambda >  \^q(\alpha) )\overset{}{\lr}1.
\end{align*}
In other words, the test described in (\ref{h0}) with statistic $\^\Lambda$ and CS bootstrap calculation of quantile is consistent.
\end{theorem}

We provide the following example under $H_1$, where $\^\theta_c$ does not converge to a constant in probability. Although we cannot get the conclusion of Theorem \ref{bootcor}, the least squared constrained statistic still converges in distribution.

\begin{example}\label{ex2}\normalfont
Let $(X_i)_{i\in \N}$ be a i.i.d. sequence such that $X_1 \overset{\dd}{=}\mathcal N (0,1)$. Define $\^\theta = \overline{ X}$, and  $H_0: \ \theta_0^2 = 1$. Clearly $H_0$ does not hold and naturally the statistic $n \ \underset{\theta^2 = 1}{\min} \|\^\theta - \theta \|^2$
goes to infinity in probability. One can find that $\^\theta_c = \text{sign}(\overline{ X})$ which does not converge. Since
\begin{align*}
\theta^*_c = \underset{\theta^2=1}{\argmin } \| \theta^*_0 - \theta\|^2  \quad \text{and}\quad \theta^*_0 = \^\theta_c + n^{-1/2}W^*,
\end{align*}
we get that $\theta^*_c= \^\theta_c$ a.s. and naturally, we do not have the 
asymptotic given by Theorem \ref{bootcor}. Besides, the convergence to a 
chi-squared distribution holds for the quantity 
$n \ \underset{\theta^2=1}{\min } \| \theta^*_0 - \theta\|^2 $. 
\end{example}

\section{Rank estimation with hypothesis testing}\label{s3}

In this section through a review of the literature about rank 
estimation, we apply the results obtained in section 
\ref{s21} to provide a consistent bootstrap procedure for the test 
described by (\ref{hyp}) associated with the statistics $\^\Lambda_1$, 
$\^\Lambda_2$ and $\^\Lambda_3$. We define $q_0=p-d_0$ the dimension of the kernel of $M_0^T$. We denote 
by $(\lambda_1,..., \lambda_{p})$ the singular values of $M_0$ arranged in descending order and we write the SVD of $M_0$ as 
\begin{align*}
M_0 = (U_1 U_0) \begin{pmatrix} D_1 & 0 \\ 0 & 0 \end{pmatrix} \begin{pmatrix} V_1^T \\V_0^T\end{pmatrix},
\end{align*} 
with $U_1\in \R^{p\times d_0}$, $U_0\in \R^{p\times q_0}$, $V_1\in \R^{H\times d_0}$, $V_0\in \R^{H\times q_0}$, and $D_1= \diag(\lambda_1,..., \lambda_{d_0})$. For $m\in \{1,\cdots,p\}$, we note $q= p -m$ and we write the SVD of $\^M$ as
\begin{align*}
\^M = (\^U_1 \^U_0) \begin{pmatrix} \^D_1 & 0 \\ 0 & \^D_0 \end{pmatrix} \begin{pmatrix} \^V_1^T \\ \^V_0^T\end{pmatrix},
\end{align*}
with $\^U_1\in \R^{p\times m}$, $\^U_0\in \R^{p\times q}$, $\^V_1\in \R^{H\times m}$, $\^V_0\in \R^{H\times q}$, $\^D_1= \diag(\^\lambda_1,..., \^\lambda_{m})$ and $\^D_0= 
\diag(\^\lambda_{m+1},..., \^\lambda_{p})$. We also introduce the orthogonal projectors
\begin{align*}
Q_1 =I-P_1= U_0 U_0^T, \  Q_2=I-P_2= V_0V_0^T,\  \^Q_1 =I-\^P_1= \^U_0 \^U_0^T \  \text{and} \  \^Q_2  =I-\^P_2=\^V_0 \^V_0^T.
\end{align*}
Whereas the link between 
$\^\Lambda_3$ and LSCE is evident, the one conecting $\^\Lambda_1$ and
$\^\Lambda_2$ to LSCE relies on the following classical lemma, whose proof is 
avoided.  
\begin{lemma}\label{lemmunproved}
Let $\^M \in \R^{p\times H} $, it holds that
\begin{align*}
\underset{\rank(M)= m}{\argmin} \| \^M-M\|_F^2 = \^P_1 \^M \^P_2, \quad \text{and} \quad \|\^M-  \^P_1 \^M \^P_2\|_F^2 = \sum_{k=m+1}^p \^\lambda_k^2 ,
\end{align*}
where $\^\lambda_1,\ldots , \^\lambda_p$ are the singular values of $\^M$ arranged in descending order, and $\^P_1$ and $\^P_2$ are orthogonal right and left singular projectors of $\^M$ associated with $\^\lambda_1,\ldots,\^\lambda_m$.
\end{lemma}

Note that in the previous lemma, $\^P_1$ and $\^P_2$ are uniquely determined if and only if $\^\lambda_{m}\neq \^\lambda_{m+1}$.

\subsection{Nonpivotal statistic}

 As stated in the introduction, the statistic $\^\Lambda_1 
=n\sum_{k=m+1}^{p} \^\lambda_k^2$ can be used to arbitrate between the 
hypotheses of (\ref{hyp}). Basically, if $H_0:\ d_0=m$ is realized, all 
the eigenvalues of the sum goes to $0$ and 
$\^\Lambda_1$ has a weighted chi-squared limiting distribution. Otherwise, at least 
one eigenvalue converges in probability to a positive number and for 
any $A>0$,  $\P(\^\Lambda_1 > A) \overset{}{\lr} 1$. The following proposition
describes the asymptotic behaviour of $\^\Lambda_1$\footnote{A similar proposition can be stated applying Proposition \ref{ssboot}. Following this way, the asymptotic depends on $g$ which is difficult to estimate for rank constraints (see Remark \ref{ex1}).}. It was stated in \cite{cook2001} and some recent extension can be found in \cite{bura2010}. Our statement goes further because we are also concerned about the estimation of the asymptotic law of $\^\Lambda_1$, i.e. the estimation of the weights that intervenes in the weighted chi-squared asymptotic law. Besides, the proof we give in the Appendix is quite simple\footnote{We no longer need the results of \cite{eaton1994} about the asymptotic behaviour of singular values.}.

\begin{proposition}\label{th1}
Under $H_0$, if (\ref{as1}) holds we have
\begin{align*}
\^\Lambda_1 \overset{d}{\lr} \sum \nu_k W_k^2
\end{align*}
where the $\nu_k$'s are the eigenvalues of the matrix $ (Q_2 \otimes 
Q_1 ) \Gamma (Q_2 \otimes Q_1 )$ and the $W_k$'s are i.i.d. standard 
Gaussian variables. If moreover (\ref{as2}) holds, we have
\begin{align*}
(\^\nu_1,..., \^\nu_{p H}) \overset{\P}{\lr} (\nu_1,..., \nu_{p H}),
\end{align*} 
where the $\^\nu_k$'s are the eigenvalues of the matrix $ (\^Q_2 \otimes \^Q_1 ) \^\Gamma (\^Q_2 \otimes \^Q_1 )$.
\end{proposition}

\begin{remark}\label{crit}\normalfont
Unlike Theorem 1 in \cite{cook2001} or Theorem 1 in 
\cite{bura2010}, we prefer to state this theorem with the quantities 
$Q_1$ and $Q_2$ rather than with $U_0$ and $V_0$. Because we do not 
assume that the kernel of $M$ has dimension $1$, the vectors that form 
$U_0$ or $V_0$ are not unique because vector spaces with dimension 
larger than $2$ have an infinite number of basis. As a consequence it does not 
make sense to estimate either $U_0$ or $V_0$. To characterize convergence of spaces, a suitable object is their associated orthogonal projectors. 
\end{remark}

In general, we do not know the asymptotic distribution of 
$\^\Lambda_1$ because it depends on $(Q_2 \otimes 
Q_1 ) \Gamma (Q_2 \otimes Q_1 )$. On the first hand, one can estimate 
consistently this matrix to get an approximation of the law of 
$\^\Lambda_1$ under $H_0$. Some conditions providing the consistency of the estimation are stated in Proposition \ref{th1}. On the other hand, one can apply the CS bootstrap to estimate the quantile of $\^\Lambda_1$ in order to test. The main advantage of such an approach is that we no longer need to have a consistent estimator of $\Gamma$ so that (\ref{as2}) is not needed anymore. Following section \ref{s21} and by using Lemma \ref{lemmunproved}, we define
\begin{align}\label{bootmatrix1}
M^*_0 = \^P_1 \^M \^P_2 + n^{-1/2} W^*\qquad \text{with} \quad   W^*|\^P \overset{\dd}{\r} W\quad \as,
\end{align}
with $W$ defined in (\ref{as1}). Accordingly, we introduce the CS bootstrap statistic 
\begin{align*}
\Lambda_1^* = n\sum_{k=m+1}^ {p} \lambda_k^{*2},
\end{align*}
with $\lambda_{m+1}^*,..., \lambda_p^*$ the smallest singular values of $M^{*}$. The following proposition is a straightforward application of Theorem \ref{bootthtest} with the submanifold $\{\rank(M)=m\}$.

\begin{proposition}\label{thboot1}
If (\ref{as1}), (\ref{bootmatrix1}) and $\^M\overset{\as}{\r} M_0$ hold, then the test described in (\ref{hyp}) with the statistic $\^\Lambda_1$ and calculation of quantile with $\^\Lambda_1^*$ is consistent.
\end{proposition}

\subsection{Wald-type statistic}
The Wald-type statistic $\^\Lambda_2 = \vecv(\^Q_1 \^M \^Q_2)^T  
[(\^Q_2\otimes \^Q_1)\^\Gamma(\^Q_2\otimes \^Q_1)]^{+} \vecv(\^Q_1 \^M \^Q_2)$ has been introduced in \cite{bura2010} to get a 
pivotal statistic\footnote{We write the expression of $\^\Lambda_2$ another way for the reasons explained in Remark \ref{crit} but one can recover the original expression by noting that for any symmetric matrix $A$, $A^+H=(AH)^+$ if $H$ is an orthonormal basis of a vector subspace of $\im(A)$.}. They obtained the following theorem for which we provide a different proof in the appendix.
 
\begin{proposition}\label{bura}
If (\ref{as1}) and (\ref{as2}) hold, we have
\begin{align*}
\^\Lambda_2 \overset{\dd}{\lr} \chi_{s}^2,
\end{align*}
with $s = \min(\rank(\Gamma),(p-d)(H-d))$.
\end{proposition}

Following (\ref{statboot}), we define the associated bootstrap statistic by
\begin{align*}
\^\Lambda_2^* = \vecv(Q_1^* M^*_0 Q_2^*)^T [(Q_2^*\otimes Q_1^*) \Gamma^{*}(Q_2^*\otimes Q_1^*) ]^+  \vecv( Q_1^* M^*_0 Q_2^*),
\end{align*}
where $M^*_0$ is defined in (\ref{bootmatrix1}), $\Gamma^*\in \R^{pH\times pH}$, $\^Q_1^*$, and $\^Q_2^*$ are the  eigenprojectors associated with the smallest eigenvalues of $M^*_0 M_0^{*T}$ and $M^{*T}_0 M^*_0$. As Proposition \ref{thboot1}, the following one is an easy application of Theorem \ref{bootthtest}.


\begin{proposition}\label{thboot2}
If  (\ref{as1}), (\ref{as2}), (\ref{bootmatrix1}), $\^M\overset{\as}{\r} M_0$ and 
$\Gamma^*\overset{\P}{\r} \Gamma$  
hold, then the test described in (\ref{hyp}) with the statistic $\^\Lambda_2$ and calculation of quantile with $\Lambda_2^*$ is consistent.
\end{proposition}

\subsection{Minimum Discrepancy approach}

Noting that $\{\rank(M)=m\}$ has co-dimension $(H-m)(p-m)$ and applying (\ref{ssboot}) we get the following proposition\footnote{See \cite{donald1997} for the original proof.}.
\begin{proposition}\label{conv3}
If (\ref{as1}), (\ref{as2}), and (\ref{as3}) hold, we have
\begin{align*}
\^\Lambda_3 \overset{\dd}{\lr} \chi_{(H-m)(p-m)}^2.
\end{align*}
\end{proposition}

In general a minimizer 
\begin{align*}
\^M_c = \underset{\rank(M) = m}{\argmin} \vecv(\^M-M)^T \^\Gamma^{-1} \vecv(\^M-M)
\end{align*}
does not have an explicit form as it was for the constrained 
matrix associated with $\^\Lambda_1$ and $\^\Lambda_2$. Therefore, we define
\begin{align}\label{bootmatrix2}
M^*_0 = \^M_c + n^{-1/2} W^*\qquad \text{with} \quad   W^*|\^P \overset{\dd}{\r} W\quad \as,
\end{align}
where $W$ is defined in (\ref{as1}). We also define the associated CS bootstrap statistic 
\begin{align*}
\Lambda_3^* = n \min_{\rank(M)=m} \vecv(M^*_0-M)^T \Gamma^{*-1} \vecv(M^*_0-M), 
\end{align*}
and applying Theorem \ref{bootthtest} we have the following result.

\begin{proposition}\label{thboot3}
If (\ref{as1}), (\ref{as2}), (\ref{as3}), (\ref{bootmatrix2}), 
$\Gamma^*\overset{\P}{\r} \Gamma$, and $\^M\overset{\as}{\r} M_0$ 
hold, then the test described in (\ref{hyp}) with the statistic 
$\^\Lambda_3$ and calculation of quantiles with $\Lambda_3^*$ is consistent.
\end{proposition}

\begin{remark}\normalfont
The set of assumptions needed to obtain Proposition \ref{conv3} is stronger than 
the ones stated in propositions \ref{th1} and \ref{bura} ensuring the convergence of $\^\Lambda_1$ and 
$\^\Lambda_2$. As a consequence this is also true for Proposition \ref{thboot3} 
with respect to propositions \ref{thboot1} and \ref{thboot2}. The main 
difference is that we add the assumption on $\Gamma$ to be non deficient. This assumption cannot be alleviated in the statement but is not as restrictive in practice. On the first hand, if $\Gamma$ is deficient the optimization under constraint has a free coordinate which implies the non-convergence of the minimizer. On the other hand, because of the semi-definite character of $\Gamma$ the projection of $\^M$ on the null space of $\Gamma$ is null. Then one can apply the proposition to the restriction of $\^M$ on the range of $\Gamma$. This is the case in the application to SDR in Section \ref{s5}.
\end{remark}

\begin{remark}\normalfont
Unlike the situation of $\^\Lambda_1$ and $\^\Lambda_2$, an optimization algorithm is 
needed to obtain $\^\Lambda_3$ and $\Lambda^*_3$, this point out an important issue of such a procedure. In \cite{cook2005}, the authors noticed that
\begin{align*}
\^\Lambda_3=n \underset{A\in H_d ,B\in \R^{d\times l}}{\min} (\vecv(\^M)- \vecv(AB) )^T \^\Gamma^{-1} ( \vecv(\^M) - \vecv(AB)) 
\end{align*}
where $H_d$ is the set of orthogonal basis lying in $\R^p$ with 
dimension $d$. We follow their algorithm in the computation of $\^\Lambda_3$ (see \cite{cook2005}, Section 3.3 for the details). 
\end{remark}

\subsection{The statistics $\^\Lambda_1$, $\^\Lambda_2$, $\^\Lambda_3$ through an example}
In the introduction, we already mentioned several drawbacks and advantages of the use of $\^\Lambda_1$, $\^\Lambda_2$, or $\^\Lambda_3$. The remark relied on both pivotality of the statistics and large matrix inversion. Here we develop another point of view related to the algebraic nature of the statistics. Facing the representation provided by Table \ref{table1}, each statistic $\^\Lambda_1$ and $\^\Lambda_2$ evaluates a different distance between $\^M$ and $\^M_c$. The first one is the distance that is optimized, but the second is another one. This has raised the issue we present here through the following example.  For the sake of clarity, we consider
\begin{align*}
\^M =\begin{pmatrix} \^\lambda_{1} &0 \\ 0 & \^\lambda_{2}  \end{pmatrix} \qquad \text{with}  \ \^\lambda_k = 
\frac 1 n \sum_{i=1}^n \lambda_{k,i},\ \text{for}\ k=1,2,\ 
\text{and}\ (\lambda_{k,i})_{k,i}\ \text{i.i.d.,}
\end{align*}
and we test $H_0: d_0=1$ against $H_1: d_0>1$. We assume that $\^\lambda_1>\^\lambda_2$, we have $\^\Lambda_1 = n \^\lambda_2^2$. Otherwise, one can show that $\^\Lambda_2= n\frac {\^\lambda_2^2 }{\^v_{2}} + o_{\P}(1) $, with $\^v_k = \overline{(\lambda_{k}-\overline{\lambda_k})^2} $. For $\^\Lambda_3$ it is clear that the minimization can be done over the diagonal matrix $\diag(\lambda_1, \lambda_2)$ and one has\begin{align*}
\^\Lambda_3 = n \ \underset{ \lambda_1 \lambda_2 = 0}{\argmin}\left\{ \frac{\^\lambda_{1}-\lambda_1}{\^v_1 } +\frac{\^\lambda_{2}-\lambda_2}{\^v_2 }    \right\}+o_{\P}(1) = n\min\left(\frac{\^\lambda_1^2}{\^v_{1}} ,\frac{\^\lambda_2^2}{\^v_{2}}\right)+o_{\P}(1).
\end{align*}
Accordingly, by Proposition \ref{thboot1}, \ref{thboot2} and \ref{thboot3}, the three tests can be summarized by 
\begin{align*}
n\^\lambda_2^2 \qquad &\text{compared to}\qquad v_{2}\chi_1^2, \\
 n\frac{\^\lambda_2^2}{\^v_{2}}  \qquad &\text{compared to}\qquad 
\chi_{2}^2,\\
 n\min\left(\frac{\^\lambda_1^2}{\^v_{1}} ,\frac{\^\lambda_2^2}{\^v_{2}}\right) \qquad &\text{compared to}\qquad 
\chi_{2}^2,
\end{align*}
where $v_k = \var(\lambda_{k,1})$. Assume there is less variance on the estimate of the smallest eigenvalue, i.e. $v_1>v_2$ such that $\frac{\^\lambda_1^2}{\^v_{1}} < \frac{\^\lambda_2^2}{\^v_{2}}$, this situation may arise when $\^\lambda_1$ and $\^\lambda_2$ have similar values but different variances. Then to conduct the test, the statistic $\frac{\^\lambda_1^2}{\^v_{1}}$ is a better choice than $\frac {\^\lambda_2^2}{\^v_{2}}$. As a consequence, unlike $\^\Lambda_1$ and $\^\Lambda_2$, the statistic $\^\Lambda_3$ appears as a coherent choice because its associated minimization takes into account the variance of the estimation.


\section{Application to sufficient dimension reduction}\label{s5}

We focus on a particularly famous method in SDR called sliced inverse regression (SIR) which has been introduced in \cite{li1991} to deal with the regression model
\begin{align}\label{cs}
Y= f(PX, \varepsilon)
\end{align}
where $\varepsilon \indep X \in \R^p$, $Y\in \R$, and $P$ is a projector on the vector space $E$ with dimension $d_0<p$, called the central subspace. The objective is to estimate $E$. If $X$ is elliptically distributed, then we have that $\Sigma^{-1}(\E[(X-\E[X])\psi(Y)]\in E $ with $\Sigma = \var(X)$, for any measurable function $\psi$. Accordingly, in order to recover the whole central subspace one needs to consider many functions $\psi$. For a given family of functions $(\psi_h)_{1\leq h\leq H}$ we define $\Psi = (\psi_1(Y),...,\psi_H(Y))^T$. Under some additional conditions \cite{portier2011}, the image of the matrix $\Sigma^{-1/2}\cov(X, \Psi(Y))$ is equal to $\Sigma^{1/2} E$. Then one can make the svd of an estimator of this matrix to obtain $d_0$ vectors that form an estimated basis of $\Sigma^{1/2}E$. Motivated by the curse of dimensionality, the estimation of $d_0$ is one of the most crucial points in SDR. To make that possible, a popular way consists in estimating the rank of $\Sigma^{-1/2} \cov(X,\Psi)$ using the hypothesis testing framework given by (\ref{hyp}) (see for example \cite{li1991}, \cite{cook2001} and \cite{cook2005}). Since we are interested in estimating the rank, we prefer to deal directly with $ \cov(X, \Psi)$ to avoid the introduction of an additional noise due to the estimation of the matrix $\Sigma$. Assume that $((X_1,Y_1),\cdots, (X_n,Y_n))$ is a i.i.d. sequence from model (\ref{cs}), denote by $\^P$ its associated empirical c.d.f. and define the quantity 
\begin{align*}
C = \E[K], \qquad \text{with}\  K = (X- \E[X])(\Psi(Y)-\E[\Psi(Y)])^T,
\end{align*}
associated with its empirical estimator
\begin{align*}
\^C = \overline{\^K},\qquad \text{with } \^K_i = (X_i-\overline{X})(\Psi_i-\overline{\Psi})^T, \quad  \text{and } \Psi_i = \Psi(Y_i).
\end{align*} 
We apply the CS bootstrap to calculate the quantiles of each statistic. Facing (\ref{bootmatrix1}) and (\ref{bootmatrix2}), we use an independent weighted bootstrap to reproduce the asymptotic law of $\sqrt n (\^C- C)$, that is we define the bootstrap matrix
\begin{align}\label{bsample}
 C^* = \^C_{c}  +  \overline{ K^*},\quad \text{with}\quad  K_i^* = w_i (\^K_i-\overline{\^K}) 
\end{align}
where $\^C_{c}$ stands for the solution of an optimization problem depending on the selected statistic $\Lambda_1$, $\Lambda_2$ or $\Lambda_3$ (see Section \ref{s3} for the details) and $(w_i)$ is a sequence of i.i.d. random variables. We also define
\begin{align*}
V=\var(\vecv(K)) \qquad \text{and} \qquad V^{*}= \frac 1 {n} \sum_{i=1}^n \vecv(K_i^*- \overline{K^*})\vecv(K_i^*- \overline{K^*})^T.
\end{align*}
To apply propositions \ref{thboot1}, \ref{thboot2}, and \ref{thboot3}, we need the following result which is of particular interest since it provides a new bootstrap procedure for SIR that is different than the one proposed in \cite{velilla2007}.

\begin{proposition}\label{bootsir}
Assume that $\E[\|X\|^2]<+\infty $, $\E[\|\Psi(Y)\|^2]$ and $\E[\|K\|_F^4]$ are finites, if moreover $(w_i)$ is a i.i.d. sequence of real random variables with mean $0$ and variance $1$, then we have
\begin{align*}
\mathcal L _{\infty}(n^{1/2}\ \overline{K^*} |\^P) = \mathcal L _{\infty} (n^{1/2}( \^C- C))\quad\text{a.s. and} \quad V^* \overset{\P}{\r} V \quad \text{conditionally a.s..}
\end{align*} 
\end{proposition}

\begin{remark}\normalfont
Taking a partition $\{I(h),\ h=1,\ldots,H\}$ of the range of $Y$ we recover the original SIR method with the family formed by the $ p_h^{-1/2} \mathds 1 _{\{Y\in I(h)\}}$'s with $p_h =\P(Y\in I(h))$. Then $C_{\text{SIR}}=\Sigma^{-1/2}\cov(X,\mathds 1) D^{-1/2}$ with $\mathds 1 =(\mathds 1 _{\{Y_i\in I(1)\}},\ldots ,\mathds 1 _{\{Y_i\in I(H)\}})^T$ and $D = \diag(p_h)$, is estimated by $\^C_{\text{SIR}}=\^\Sigma^{-1/2}\overline{(X-\overline{X}) \mathds 1^T}\^D^{-1/2}$ with $\^D = \diag (\^p_h)$, $\^p_h =\overline{\mathds 1 _{\{Y\in I(h)\}}}$, $\^\Sigma = \overline{(X-\overline{X})(X-\overline{X})^T}$. We have the expansion
\begin{multline*}
n^{-1/2}(\^C_{\text{SIR} }- C_{\text{SIR}}) =  n^{-1/2} \Sigma^{-1/2} (\overline{(X-\E[X])  \mathds 1 ^T}- \cov(X,\mathds 1 ) ) D^{-1/2}\\ -\Sigma^{-1/2}n^{-1/2} (\^\Sigma^{1/2}-\Sigma^{1/2})C_{\text{SIR}} -  C_{\text{SIR}}n^{-1/2} (\^D_p^{1/2}-D_p^{1/2})D_p^{-1/2}+ o_{\P}(1).
\end{multline*}
As a consequence, the matrix $\Sigma^{-1/2}$ and the weights $p_h$'s are playing an important role on the asymptotic of the matrix SIR. They introduce some other terms in the asymptotic distribution and clearly the simple bootstrap presented before does not work for SIR as it was originally defined. Even if we believe that a more evolved weighted bootstrap works to bootstrap $\sqrt n (\^C_{\text{SIR} }- C_{\text{SIR}}) $, we emphasize that it may be less accurate than the one we propose since it complicates the asymptotic without being necessary for testing the rank. 
\end{remark}

Recall that $m$ is a non-negative integer, for $k\in \{1,2,3\}$ and $B\in \N^*$ we calculate independent copies $\Lambda_{k,1}^{*},...,\Lambda_{k,B}^{*}$ with the CS bootstrap algorithm corresponding to each statistic. Then we estimate the quantile with 
\begin{align*}
q_k^*(\alpha)=  \inf_{t\in \R}\{ F_k^* (t) >\alpha \} = \Lambda_{k,(\lceil B\alpha\rceil )}^*,\qquad \text{where} \ F_k^*(t) = \frac 1 B \sum_{b=1}^B \mathds 1 _{\{\Lambda_{k,b}^* \leq t\}},
\end{align*}
$\lceil \cdot \rceil$ is the integer ceiling function and $\Lambda_{k,(\cdot)}^*$ stands for the rank statistic associated to the sample $\Lambda_{k,1}^*\ldots\Lambda_{k,B}^*$. On the first hand, we conduct the test described by (\ref{hyp}) using the CS bootstrap, i.e.  
\begin{align}\label{test}
H_0\text{ is rejected if} \qquad \^\Lambda_k > \^q_k^*(\alpha).  
\end{align}
On the other hand, the traditional test is conducted by comparing the statistic $\^\Lambda_2$ and $\^\Lambda_3$ to the quantile of their asymptotic law respectively given by propositions \ref{bura} and \ref{conv3}. For $\^\Lambda_1$, in general the limit in law is quite complicated\footnote{When the predictors are normally distributed, it has been shown that $\^\Lambda_1$ is asymptotically chi-squared distributed (see \cite{cook2001}). The authors also pointed out that it was less robust than the weighted chi-squared asymptotic as soon as the predictors distribution deviates from normality. As a result, we keep in the nonparametric framework by avoiding such asymptotic in this simulation study.} (see Proposition \ref{th1}), so that we use approximations: the Wood's approximation (see \cite{wood1989}) as it is computed in the R software, an adjusted version $\^\Lambda_{1,\text{adj.}} = \^\Lambda_1 / a \overset{\dd}{\r} \chi_{b}^2$, with $a= \sum_{k=1}^{s}\omega^2 /\sum_{k=1}^{s} \omega_k $, $b= (\sum_{k=1}^{s}\omega_k)^2 /\sum_{k=1}^{s} \omega_k^2 $, and a re-scaled version $\^\Lambda_{1,\text{sc}} =\^\Lambda_1 /c \overset{\dd}{\r} \chi_s^2$, $c= \overline{\omega}$ (see \cite{bentler2000} for these two corrections). 

In all the simulations we compute the matrix $\^C$ by taking $\Psi(t) = (\mathds 1 _{\{y\in I(1),\ldots , y\in I(H)\}})$ where the $I(h)$'s form an equi-partition of the range of the data $Y_1,\ldots , Y_n$. In the whole study we put $(p,H)= (6,5)$, $B=1000$ and we consider $n=50,100,200,500$. Although the parameter $H$ does not really affect the SIR method, we choose it globally good with respect to all the situations.

\begin{figure}\centering
\includegraphics[width=7.5cm]{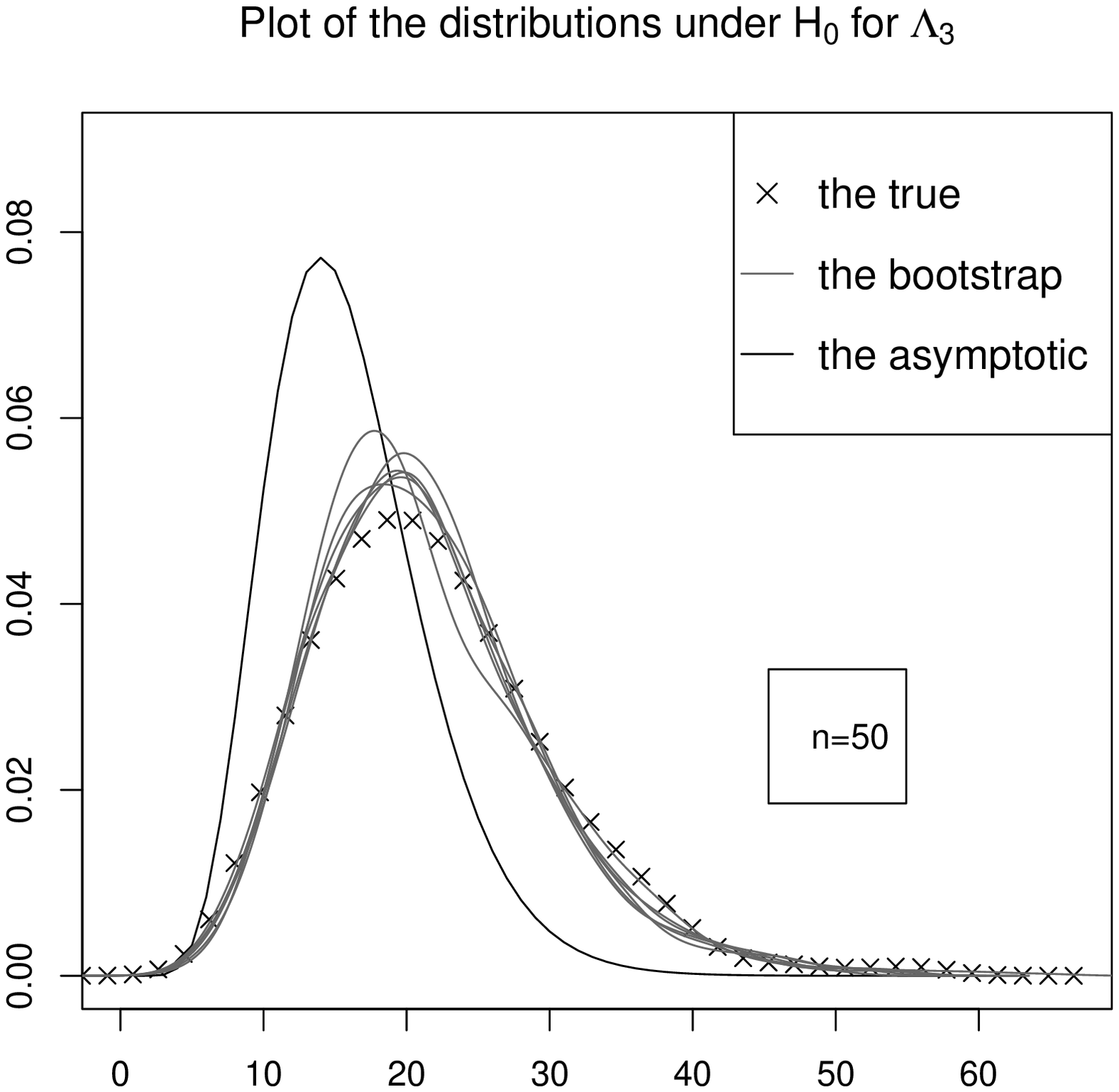} 
\includegraphics[width=7.5cm]{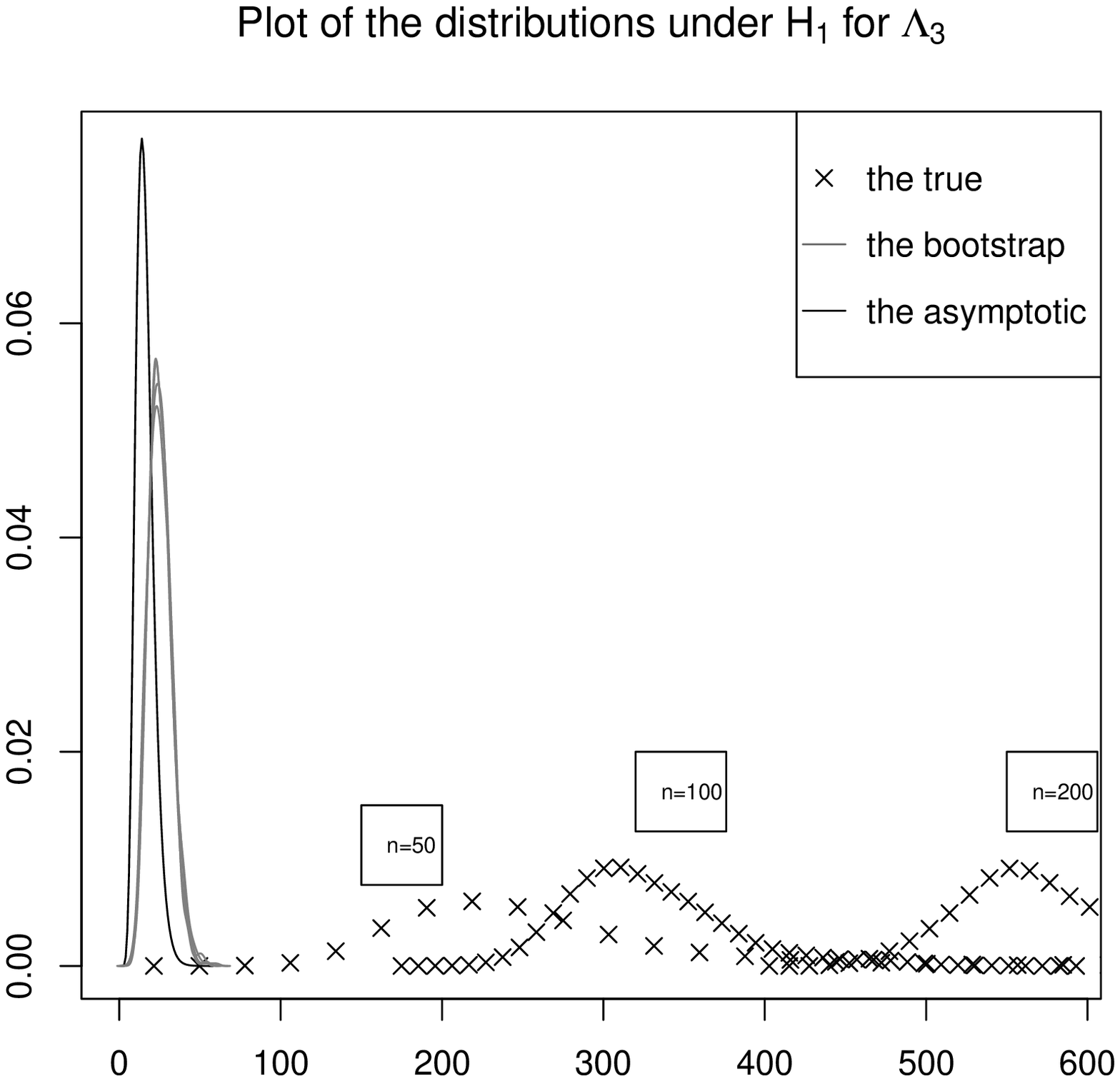}
\caption{Plot of the asymptotic distribution, and the estimated distribution of the statistic and the bootstrap statistic for $\^\Lambda_3$ in the case of Model \cRM{1}.}\label{fig1}
\end{figure}

The first model we study is the following standard model:
\begin{align*}
\text{Model \cRM{1}:} \qquad\qquad Y= X_1+.1e\qquad \qquad\text{with } e\indep X,\quad X\overset{\text{d}}{=}\mathcal N (0,I), \quad e\overset{\text{d}}{=}\mathcal N (0,1).
\end{align*}
In order to highlight guidelines (\ref{gui1}) and (\ref{gui2}), we produce in figure \ref{fig1} two graphics each representing situation under $H_1$ and $H_0$ for the statistic $\^\Lambda_3$. Similar graphics dealing with $\^\Lambda_2$ have been drawn but are not presented here. 
On the first one we see that even if the sample is under $H_1$ the bootstrap distribution reflects $H_0$. As a consequence, guideline (\ref{gui1}) is satisfied and the power of the bootstrap test is going to $1$. The second graph shows that the statistic distribution is closer to the bootstrap distribution than its asymptotic distribution. This has no reason to occur when the statistic is not pivotal (see the introduction and \cite{hall1992} for the details). As a consequence, we believe that this good fitting is due to Guideline \ref{gui2}.

\begin{figure}\centering
\includegraphics[width=8.5cm]{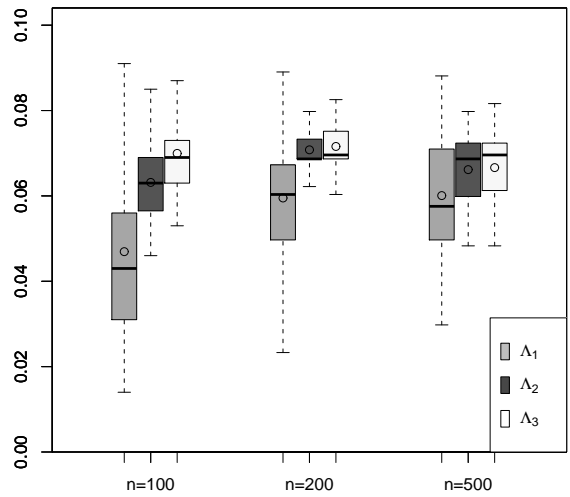} 
\caption{Bowplot over $100$ samples of $\^q(\alpha)$ for $\^\Lambda_1$, $\^\Lambda_2$, $\^\Lambda_3$ and $\alpha=0.95$ in the case of Model \cRM{1} for different values of $n$ .}\label{fig2}
\end{figure}

In figure \ref{fig2} we analyse the asymptotic distribution of $\^q(\alpha)$ in model \cRM{1} for each statistic. To measure the error we consider the behaviour of 
\begin{align*}
F_n (\^q(\alpha)),
\end{align*}
which is optimally equal to $1-\alpha$. To make that possible, $F_n$ is estimated with a large sample size so that the estimation error is negligible. Then we run over $100$ samples the CS bootstrap to provide, for each sample, a bootstrap estimation of the quantile $\^q(\alpha)$. The associated boxplot for $n = 100$, $200$, $500$ are provided in Figure \ref{fig2}. As a consequence, we may notice that the behaviour
of $\^\Lambda_2$ and $\^\Lambda_3$ are quite similar facing the one of $\^\Lambda_1$. Even if every boxplot argues for convergence to $1-\alpha$, testing with $\^\Lambda_1$ seems a better choice when n is small because of a quasi immediate convergence of the bias. When n increase, this is no longer evident because the variance of either $\^\Lambda_2^*$ or$\^\Lambda_3^*$ is smaller.

Furthermore, we go into details in Table \ref{tab1} by running Model \cRM{1} over $5000$ samples. For each of them and every statistic, we conduct the bootstrap test (\ref{test}) and its traditional version. The table presents for each $m\leq d_0$, the proportion of rejected tests. This corresponds to either estimate of the power or estimate of the level.

\newcolumntype{A}{ >{\columncolor[gray]{.8}}c }

\begin{table}\centering \small
\begin{tabular}{l  c | c c c A || c A c A}
 \hline
    \multirow{2}{*}{$n$} & \multirow{2}{*}{$m$}  & \multicolumn{4}{ >{\columncolor[gray]{.7}}c }{$\^\Lambda_1$}& \multicolumn{2}{ >{\columncolor[gray]{.7}}c }{$\^\Lambda_2$} &\multicolumn{2}{ >{\columncolor[gray]{.7}}c }{$\^\Lambda_3$} \\ 
 & & \multicolumn{1}{ >{\columncolor[gray]{.7}}c }{Wood} & \multicolumn{1}{ >{\columncolor[gray]{.7}}c }{Resc.} &\multicolumn{1}{ >{\columncolor[gray]{.7}}c }{Adj.}& \multicolumn{1}{ >{\columncolor[gray]{.7}}c }{CB $\^\Lambda_1$}& \multicolumn{1}{ >{\columncolor[gray]{.7}}c }{$\^\Lambda_2$}&\multicolumn{1}{ >{\columncolor[gray]{.7}}c }{CB $\^\Lambda_2$}&\multicolumn{1}{ >{\columncolor[gray]{.7}}c }{$\^\Lambda_3$} &\multicolumn{1}{ >{\columncolor[gray]{.7}}c }{CB $\^\Lambda_3 $}\\ \hline
  \hline
\multirow{2}{*}{$50$}& $0$ & 0.9988& 0.9998& 0.9988& 0.9988  &1.0000 &1.0000  &1.0000& 1.0000 \\ 
& $1$ & 0.0326 &0.0590 &0.0336 &0.0494  &0.3466 & 0.0744&0.3098  &0.07\\  \hline
\multirow{2}{*}{$100$}& $0$ & 1.0000& 1.0000 &1.0000 &1.0000 &1.0000 &1.0000 &1.0000&1.0000\\ 
& $1$ &0.0386 &0.052 &0.0388& 0.0456& 0.1494 &  0.0676&0.1466 &0.0722\\  \hline
\multirow{2}{*}{$200$}& $0$ & 1.0000& 1.0000 &1.0000 &1.0000 &1.0000 &1.0000&1.0000&1.0000 \\
& $1$ &0.0474 &0.055 &0.0476& 0.0514 &0.096&  0.0646&0.0954 &0.0664\\  \hline
\multirow{2}{*}{$500$}& $0$ & 1.0000& 1.0000 &1.0000 &1.0000 &1.0000 &1.0000&1.0000&1.0000\\ 
& $1$ &0.0492 &0.0514 &0.0494& 0.0516 &0.0656 &0.0584&0.0654 & 0.0584\\  \hline
\end{tabular}\caption{Estimated levels and power in Model \cRm{1} for $\alpha=5\%$.}\label{tab1}
\end{table}

Although it has not the best power, the clear winner is the tests based on $\^\Lambda_1$. Inside this group, for any sample number, the bootstrap and the rescaled version are the closest to the nominal level. Concerning $\^\Lambda_2$ and $\^\Lambda_3$ the result are quite impressive when $n$ is small: for $n=100$, whereas traditional testing makes a type \cRM{1} error $30$\% of the time, the bootstrap testing goes wrong around $7$\%. This confirms observation on the second graph of Figure \ref{fig1}.

\begin{table}\centering\small
\begin{tabular}{l  c | c c c A || c A c A}
 \hline
    \multirow{2}{*}{$n$} & \multirow{2}{*}{$m$}  & \multicolumn{4}{ >{\columncolor[gray]{.7}}c }{$\^\Lambda_1$}& \multicolumn{2}{ >{\columncolor[gray]{.7}}c }{$\^\Lambda_2$} &\multicolumn{2}{ >{\columncolor[gray]{.7}}c }{$\^\Lambda_3$} \\ 
 & & \multicolumn{1}{ >{\columncolor[gray]{.7}}c }{Wood} & \multicolumn{1}{ >{\columncolor[gray]{.7}}c }{Resc.} &\multicolumn{1}{ >{\columncolor[gray]{.7}}c }{Adj.}& \multicolumn{1}{ >{\columncolor[gray]{.7}}c }{CB $\^\Lambda_1$}& \multicolumn{1}{ >{\columncolor[gray]{.7}}c }{$\^\Lambda_2$}&\multicolumn{1}{ >{\columncolor[gray]{.7}}c }{CB $\^\Lambda_2$}&\multicolumn{1}{ >{\columncolor[gray]{.7}}c }{$\^\Lambda_3$} &\multicolumn{1}{ >{\columncolor[gray]{.7}}c }{CB $\^\Lambda_3 $}\\ \hline
  \hline
\multirow{2}{*}{$50$}& $0$ &  0.9646 &0.9928& 0.9656 &0.9682 &1.0000& 1.0000 & 1.0000& 1.0000 \\ 
& $1$ & 0.0318 &0.0628 &0.0324 &0.0496& 0.3412& 0.0588 & 0.3042 &0.0628\\  \hline
\multirow{2}{*}{$100$}& $0$ &  0.9996& 1.0000 &0.9996& 0.9996& 1.0000 &1.0000  &1.0000 &1.0000\\ 
& $1$ & 0.0336& 0.0486 &0.0344 &0.0412& 0.1516 &0.0696 & 0.1432 &0.0718\\  \hline
\multirow{2}{*}{$200$}& $0$ & 1.0000& 1.0000 &1.0000 &1.0000 &1.0000 &1.0000&1.0000&1.0000 \\
& $1$ & 0.0378& 0.0486& 0.038& 0.0424 &0.0844 &0.0602 &0.0832& 0.0604\\  \hline
\multirow{2}{*}{$500$}& $0$ & 1.0000& 1.0000 &1.0000 &1.0000 &1.0000 &1.0000&1.0000&1.0000\\ 
& $1$ & 0.0454& 0.0502 &0.0458& 0.0474& 0.0638 &0.0606& 0.0634 &0.0608\\  \hline
\end{tabular}\caption{Estimated levels and power in Model \cRm{1}a for $\alpha=5\%$.}\label{tab2}
\end{table}

\begin{table}\centering\small
\begin{tabular}{l  c | c c c A || c A c A}
 \hline
    \multirow{2}{*}{$n$} & \multirow{2}{*}{$m$}  & \multicolumn{4}{ >{\columncolor[gray]{.7}}c }{$\^\Lambda_1$}& \multicolumn{2}{ >{\columncolor[gray]{.7}}c }{$\^\Lambda_2$} &\multicolumn{2}{ >{\columncolor[gray]{.7}}c }{$\^\Lambda_3$} \\ 
 & & \multicolumn{1}{ >{\columncolor[gray]{.7}}c }{Wood} & \multicolumn{1}{ >{\columncolor[gray]{.7}}c }{Resc.} &\multicolumn{1}{ >{\columncolor[gray]{.7}}c }{Adj.}& \multicolumn{1}{ >{\columncolor[gray]{.7}}c }{CB $\^\Lambda_1$}& \multicolumn{1}{ >{\columncolor[gray]{.7}}c }{$\^\Lambda_2$}&\multicolumn{1}{ >{\columncolor[gray]{.7}}c }{CB $\^\Lambda_2$}&\multicolumn{1}{ >{\columncolor[gray]{.7}}c }{$\^\Lambda_3$} &\multicolumn{1}{ >{\columncolor[gray]{.7}}c }{CB $\^\Lambda_3 $}\\ \hline
  \hline
\multirow{2}{*}{$50$}& $0$ & 1.0000& 1.0000 &1.0000 &1.0000 &1.0000 &1.0000 &1.0000&1.0000 \\ 
& $1$ &0.034 &0.1072 &0.034& 0.0378 &0.2122 &0.0396 &0.1394 &  0.015\\  \hline
\multirow{2}{*}{$100$}& $0$ & 1.0000& 1.0000 &1.0000 &1.0000 &1.0000 &1.0000 &1.0000&1.0000\\ 
& $1$ &0.037 &0.0904 &0.0374 &0.0404 &0.0986 &0.0572 &0.0614 & 0.0284\\  \hline
\multirow{2}{*}{$200$}& $0$ & 1.0000& 1.0000 &1.0000 &1.0000 &1.0000 &1.0000&1.0000&1.0000 \\
& $1$ &0.0484& 0.096& 0.0488 &0.0518& 0.0708 &0.066 &0.056 & 0.0506\\  \hline
\multirow{2}{*}{$500$}& $0$ & 1.0000& 1.0000 &1.0000 &1.0000 &1.0000 &1.0000&1.0000&1.0000\\ 
& $1$ &0.0486 &0.0912& 0.0486& 0.0490 &0.0598 &0.0664 &0.0612  & 0.0674\\  \hline
\end{tabular}\caption{Estimated levels and power in Model \cRm{1}b for $\alpha=5\%$.}\label{tab3}
\end{table}

In Table \ref{tab2} and Table \ref{tab3} we consider the same model than Model \cRM{1} excepted that we change the distribution of the predictors: in Model \cRM{1}a, $X$ has independent coordinates with a student distribution with $5$ degrees of freedom, in Model \cRM{1}b, $X \overset{\text{d}}{=} .1 X_1 \epsilon + X_2 (1-\epsilon)$ with $\epsilon \overset{\text{d}}{=}\mathcal B (1/2) $, $X_1\overset{\text{d}}{=}\mathcal N( (6,0,\cdots,0),I) $, $X_2\overset{\text{d}}{=}\mathcal N(0,I) $. For this two models, we have similar conclusions than model \cRM{1} with two new things. First, the rescaled version is not robust to the distribution of the predictors (Table \ref{tab3}). Second, the algorithm employed to optimized $\^\Lambda_3$ could failed at very small sample size.

We introduce a non linear relationship by considering the model
\begin{align*}
\text{Model \cRM{2}:}\qquad \qquad Y= \text{tanh} (X_1)+.1e \qquad \qquad  \text{with } e\indep X,\quad X\overset{\text{d}}{=}\mathcal N (0,I), \quad e\overset{\text{d}}{=}\mathcal N (0,1).
\end{align*}
In Table \ref{tab4}, we present similar results as in tables \ref{tab2}-\ref{tab4} with the difference that the nominal level is
$\alpha = 1\%$ in order to highlight differences in the power of each test. Again, the CS bootstrap induces a large improvement of the accuracy of the test with $\^\Lambda_2$ and $\^\Lambda_3$. At $n = 50$, the test
based on $\^\Lambda_1$ is less powerful than the others but it is more accurate under $H_0$. The winner remains the CS bootstrap with $\^\Lambda_1$. A new important things is that at $n = 500$, it seems better to use the CS bootstrap with $\^\Lambda_2$ and $\^\Lambda_3$. Actually this is due to the variance of the formers which is smaller than the variance of $\Lambda_1^*$ as it was already highlighted in Figure \ref{fig2}.

\begin{table}\centering\small
\begin{tabular}{l  c | c c c A || c A c A}
 \hline
    \multirow{2}{*}{$n$} & \multirow{2}{*}{$m$}  & \multicolumn{4}{ >{\columncolor[gray]{.7}}c }{$\^\Lambda_1$}& \multicolumn{2}{ >{\columncolor[gray]{.7}}c }{$\^\Lambda_2$} &\multicolumn{2}{ >{\columncolor[gray]{.7}}c }{$\^\Lambda_3$} \\ 
 & & \multicolumn{1}{ >{\columncolor[gray]{.7}}c }{Wood} & \multicolumn{1}{ >{\columncolor[gray]{.7}}c }{Resc.} &\multicolumn{1}{ >{\columncolor[gray]{.7}}c }{Adj.}& \multicolumn{1}{ >{\columncolor[gray]{.7}}c }{CB $\^\Lambda_1$}& \multicolumn{1}{ >{\columncolor[gray]{.7}}c }{$\^\Lambda_2$}&\multicolumn{1}{ >{\columncolor[gray]{.7}}c }{CB $\^\Lambda_2$}&\multicolumn{1}{ >{\columncolor[gray]{.7}}c }{$\^\Lambda_3$} &\multicolumn{1}{ >{\columncolor[gray]{.7}}c }{CB $\^\Lambda_3 $}\\ \hline
  \hline
\multirow{2}{*}{$50$}& $0$ &0.9308& 0.9884& 0.9428&0.9448 &1.0000&&0.9988 1.0000   &0.9988 \\ 
& $1$ & 0.0036 &0.0148 &0.0050&0.0086& 0.1816 &0.0148& 0.1404  &0.0130\\  \hline
\multirow{2}{*}{$100$}& $0$ & 1.0000& 1.0000 &1.0000 &1.0000 &1.0000 &1.0000 &1.0000&1.0000\\ 
& $1$ & 0.0072 &0.0122 &0.0082 &0.0096 & 0.0536&  0.02& 0.0496 &0.021\\  \hline
\multirow{2}{*}{$200$}& $0$ & 1.0000& 1.0000 &1.0000 &1.0000 &1.0000 &1.0000&1.0000&1.0000 \\
& $1$ & 0.0076& 0.0114& 0.0086&0.0102 & 0.0252 &0.0192 &0.0248   &0.02\\  \hline
\multirow{2}{*}{$500$}& $0$ & 1.0000& 1.0000 &1.0000 &1.0000 &1.0000 &1.0000&1.0000&1.0000\\ 
& $1$ & 0.0068& 0.0076 &0.007& 0.0082 & 0.012  &0.011&0.012 &0.011\\  \hline
\end{tabular}\caption{Estimated levels and power in Model \cRm{2} for $\alpha=1\%$.}\label{tab4}
\end{table}

We conclude by increasing difficulty considering the following model, introduced in \cite{li1991},
\begin{align*}
\text{Model \cRM{3}:} \qquad Y= \frac {X_1}{.5+(X_2+2)^2}+e\qquad \qquad e\indep X,\qquad X\overset{\text{d}}{=} \mathcal N (0,I) 
\end{align*}
We still present in Table \ref{tab5} the estimated level and power with the nominal level $\alpha= 2 \%$ for each test. For such a model the conclusions are quite mitigated because it induces a trade-off between high power and accurate level. Indeed when $n$ is small, the better powers are provided by the
traditional tests with $\^\Lambda_2$ and $\^\Lambda_3$. Nevertheless the more accurate levels can be found looking at the CS bootstrap with $\^\Lambda_2$ (n = 100) or $\^\Lambda_1$ (n = 200). Moreover the tests associated to $\^\Lambda_1$
without bootstrap are the worst concerning this model.
Accordingly, the simulation study highlighted the good behaviour of the CS bootstrap: in
every model it improves the accuracy of the traditional test for each statistic. One may remember that the bias of the CS bootstrap with $\^\Lambda_1$ has the faster rate of convergence with respect to the
CS bootstrap of $\^\Lambda_2$ or $\^\Lambda_3$. Otherwise, the variance of $\^\Lambda_1^*$
may be greater than the variance of $\^\Lambda_2^*$ or $\^\Lambda_3^*$. Finally, for the simple models it seems better to use the CS bootstrap with the statistic $\^\Lambda_1$.

\begin{table}\centering\small
\begin{tabular}{l  c | c c c A || c A c A}
\\ 
 \hline
    \multirow{2}{*}{$n$} & \multirow{2}{*}{$m$}  & \multicolumn{4}{ >{\columncolor[gray]{.7}}c }{$\^\Lambda_1$}& \multicolumn{2}{ >{\columncolor[gray]{.7}}c }{$\^\Lambda_2$} &\multicolumn{2}{ >{\columncolor[gray]{.7}}c }{$\^\Lambda_3$} \\ 
 & & \multicolumn{1}{ >{\columncolor[gray]{.7}}c }{Wood} & \multicolumn{1}{ >{\columncolor[gray]{.7}}c }{Resc.} &\multicolumn{1}{ >{\columncolor[gray]{.7}}c }{Adj.}& \multicolumn{1}{ >{\columncolor[gray]{.7}}c }{CB $\^\Lambda_1$}& \multicolumn{1}{ >{\columncolor[gray]{.7}}c }{$\^\Lambda_2$}&\multicolumn{1}{ >{\columncolor[gray]{.7}}c }{CB $\^\Lambda_2$}&\multicolumn{1}{ >{\columncolor[gray]{.7}}c }{$\^\Lambda_3$} &\multicolumn{1}{ >{\columncolor[gray]{.7}}c }{CB $\^\Lambda_3 $}\\ \hline
  \hline
\multirow{3}{*}{$50$}& $0$ &  0.9950 &0.9992 &0.9962& 0.9960 &1.0000& 0.9966 &1.0000 &0.9966 \\ 
& $1$ & 0.3750 &0.5342 &0.3990& 0.4676 &0.9074 &0.5066& 0.8344& 0.3270\\
&$2$ &   0.0078 &0.0156 &0.0086& 0.0240& 0.0620 &0.0164& 0.0344& 0.0136\\  \hline
\multirow{3}{*}{$100$}& $0$ &1.0000 &1.0000& 1.0000 &1.0000& 1.0000 &1.0000 &1.0000 &1.0000 \\ 
& $1$ & 0.9330& 0.9556 &0.9368& 0.9446 &0.9952& 0.9842& 0.9934& 0.9806\\
&$2$ &   0.0134 &0.0176& 0.0138& 0.0210& 0.0306& 0.0228& 0.0266& 0.0278\\  \hline
\multirow{3}{*}{$200$}& $0$ & 1.000& 1.0000 &1.000 &1.0000& 1.0000 &1.0000 &1.0000& 1.0000 \\ 
& $1$ &1.000 &1.0000 &1.000 &1.0000 &1.0000 &1.0000 &1.0000 &1.0000\\
&$2$ &  0.0154 &0.0182& 0.0158 &0.0198 &0.025 &0.024 &0.0244 &0.026\\  \hline
\multirow{3}{*}{$500$}& $0$ &1.0000 &1.000& 1.0000 &1.0000 &1.0000 &1.000 &1.0000& 1.0000\\
&$1$&  1.0000 &1.000& 1.0000 &1.0000& 1.0000 &1.0000 &1.0000 &1.0000\\
&$2$ &  0.0184 &0.0194& 0.0184& 0.02 &0.0228 &0.0228 &0.0228& 0.023\\ \hline
\end{tabular}\caption{Estimated levels and power in Model \cRm{2} for $\alpha=2\%$.}\label{tab5}
\end{table}

 \section{Concluding remarks}
Along this study, we found that the main advantages of the CS bootstrap are:
\begin{enumerate}
\item Alternative to the asymptotic comparison. This argument is even stronger since the asymptotic
law can be unknown (or difficult to estimate) or the asymptotic law remains too much
different from the statistic law (e.g. large matrix inversion).
\item By Theorem 4, which provides its consistency, the CS bootstrap works under mild assumptions.
Essentially, we ask the manifold to be locally smooth, and we require a bootstrap
18
of the unconstrained estimator.
\item The CS bootstrap is computationally as simple than the considered statistic.
\item In the case of rank testing, the CS bootstrap clearly improves the accuracy of traditional
testing (cf. the simulation study).
\end{enumerate}

Besides, there exists some natural extensions of the previous work. First although it is
suitable for testing, the form of the objective function bQ is quiet restrictive. For example, we
believe that the CS bootstrap could be extended to M and Z estimation. Secondly, conditions
that guarantee
\begin{align*}
\^q(\alpha)=q_n(\alpha) +o_{\P}(n^{-1/2})
\end{align*}
have not been provided yet. This would valid theoretically the use of the
CS bootstrap with respect to traditional testing.

\section*{Appendix}

\subsection*{Proof of Lemma \ref{bootth}}
The whole proof is made conditionally on the sample. By definition of 
$\^\theta_c$, with high probability, $A^*$ is full rank 
for $n$ large enough, we have
\begin{align}\label{ineg} 
 \| A^{*1/2} (\theta_c^*- \theta_c)\|\leq  \| A^{*1/2} 
 (\theta^*_c-\theta^*_0)\|+\| A^{*1/2} (\theta^*_0- \theta_c)\|\leq 2\| A^{*1/2} (\theta^*_0- \theta_c)\|.
\end{align}
Then since $\theta^*_0 -\^\theta_c\overset{\P}{\r} 0 $, $\^\theta_c\overset{}{\r}\theta_c$ and because $A^* \overset{\P}{\r} A$ is full rank, one gets that 
$\theta^*_c\overset{\P}{\r} \theta_c$. Therefore, since $\theta_c$ is 
$\mathcal M$-nonsingular and reffering to Definition \ref{const}, we get
\begin{align*}
\underset{\theta\in \mathcal M}{\argmin}\ \| 
\Gamma^{*1/2}(\theta^*_0-\theta)\| = \underset{g(\theta)=0}{\argmin}\ \| 
\Gamma^{*1/2}(\theta^*_0-\theta)\|,
\end{align*}
with $g$ continuously differentiable on $\theta_c$ and $J_g(\theta_c)$ full rank. By assumption on $g$, $\theta_c^*$, at least for $n$ large enough, satisfies the first order conditions, that are
 \begin{align*}
\left\{ \begin{array}{c}  A^*(\theta^*_0-\theta_c^*) - J_g^T (\theta_c^*) \lambda^*_n =0 \\g(\theta_c^*)=0 \end{array}\right.  
\end{align*}
where $\lambda_n^*$ is the Lagrange multiplier. Using a Taylor expansion of $g$ around $\^\theta_c$, we get $g(\theta_c^*)= g(\^\theta_c) + J_g^T(\^\theta_c)(\theta_c^*-\^\theta_c)+ o_{\P}(\|\theta_c^*-\^\theta_c\|)$, and with the previous equations we have
\begin{align*}
 \begin{pmatrix} A^* & J_g^T(\theta_c^*)\\J_g (\^\theta_c) &0 \end{pmatrix} \begin{pmatrix} \theta_c^*-\^\theta_c \\  {\lambda_n^*}  \end{pmatrix} = \begin{pmatrix} A^*(\theta^*_0-\^\theta_c)\\ o_{\P}(\|\theta_c^*-\^\theta_c\|) \end{pmatrix}.
\end{align*}
Now by Slutsky's lemma, we get
\begin{align*}
\begin{pmatrix} A & J_g^T(\theta_c)\\J_g (\theta_c) &0 \end{pmatrix} \begin{pmatrix} n^{1/2} ( \theta_c^*-\^\theta_c) \\ n^{1/2}
{\lambda_n^*} \end{pmatrix}  = n^{1/2} \begin{pmatrix} 
A(\theta^*_0-\^\theta_c)\\ 0 \end{pmatrix}+ o_{\P}(1),
\end{align*}
and the conclusion follows by multiplying on the left by the matrix 
\begin{align*}
\begin{pmatrix}A^{-1} - P A^{-1}, & A^{-1} J^T_g(\theta_c) 
(J_g(\theta_c)A^{-1} J^T_g(\theta_c))^{-1} \end{pmatrix}
\end{align*}
with $P= A^{-1} J^T_g(\theta_c)(J_g(\theta_c) A^{-1} J^T_g(\theta_c))^{-1} J_g(\theta_c)$.
 
\qed

\subsection*{Proof of Theorem \ref{bootthtest} }

The proof is divided in two parts each corresponding to the level and 
the power of the test. Assume $H_0$ and define $F_n$ and $F_{\infty}$ respectively as the c.d.f. of $\^\Lambda$ 
and the weak limit of $F_n$. Note that we can apply 
Proposition \ref{ssboot2} to get
\begin{align*}
n^{1/2} \vecdeux{\^\theta - \theta_0}{\^\theta_c - \theta_0} = 
n^{1/2} \vecdeux{I}{I-P} (\^\theta - \theta_0) +o_{\P}(1),
\end{align*}
and Theorem \ref{bootcor} to get conditionally a.s. 
\begin{align*}
n^{1/2}\vecdeux{\theta^*_0 - \^\theta_c}{\theta^*_c - \^\theta_c} = 
n^{1/2} \vecdeux{I}{I-P} (\theta^*_0 - \^\theta_c) +o_{\P}(1).
\end{align*}
with $P$ detailed in the statement of Proposition \ref{ssboot2}. Using (\ref{minchi2}), (\ref{statboot}) and Slutsky's theorem we have
\begin{align*}
\mathcal L_{\infty} (\Lambda^* |\^P) = \mathcal L_{\infty} 
(\^\Lambda)\quad \as.
\end{align*}
In other words, with probability $1$, $\^F$ converges pointwise to 
$F_{\infty}$. As in \cite{vandervaart1998} chapter 23, Lemma 3, consider $\Delta$ the set of discontinuity of $F_{\infty}^{-1}$. For every $\alpha \in ]0,1[ \backslash \Delta$, we have 
$\^q(\alpha) \overset{} {\lr}  q(\alpha)$ a.s. (see for instance \cite{vandervaart1998}, chapter 21). Using Slutsky's theorem, we get $\mathcal L_{\infty} (\^\Lambda-\^q(\alpha) ) = \mathcal L_{\infty} (\^\Lambda-q(\alpha))$, accordingly
\begin{align*}
\P(\^\Lambda \leq \^q(\alpha) )\overset{}{\lr} F_{\infty} (q(\alpha))\quad \text{for all }\alpha\in ]0,1[\backslash \Delta.
\end{align*}
Because $F_{\infty}$ is continuous $F_{\infty}(q(\alpha))= \alpha$. Since $F_{\infty}$ is non-decreasing, $\Delta$ is denumerable, since $\alpha \mapsto  \P(\^\Lambda \leq \^q(\alpha) )$ is non-decreasing with continuous limit, the convergence is uniform and so
holds for every $\alpha\in ]0,1[$. This concludes the proof for 
the level. It remains to show that the power of the test goes to $1$. Assume $H_1$ and
let $\alpha\in ]0,1[$, the statistic $\^\Lambda$ goes to infinity in 
probability and it suffices to show that with probability $1$ the bootstrap quantile $\^q(\alpha)$ remains bounded. This means exactly that conditionally a.s. the sequence $\Lambda^*$ is tight. Note that conditionally a.s. we have
\begin{align*}
\Lambda^* \leq n\| A^{*1/2} (\^\theta_c - \theta^*_0)\|^2 = \widetilde \Lambda^*, 
\end{align*}
where $\widetilde \Lambda^*$ converges in distribution by (\ref{bootcond1}), and is therefore tight.
\qed

\subsection*{Proof of Proposition \ref{th1}}
We have
\begin{align}\label{héhé}
\^\Lambda_1 = \| n^{1/2} \^Q_1 \^M \^Q_2\|^2_F=\|n^{1/2}\vecv(\^Q_1 \^M \^Q_2)\|^2.
\end{align}
By the Delta method and because $H_0$ is realized, we can apply 
convergence results about eigenprojectors to both matrices $\^M^T \^M$ and $\^M \^M^T$ to 
obtain the $\sqrt n$-convergence for $\^Q_1$ and $\^Q_2$. Then we write
\begin{align*}
n^{1/2} \^Q_1 \^M \^Q_2 &= n^{1/2}\^Q_1 (\^M-M) \^Q_2 + n^{1/2}(\^Q_1-Q_1) M (\^Q_2-Q_2)\\
&=  n^{1/2} Q_1 (\^M-M) Q_2 + O_{\P} (n^{-1/2}),
\end{align*}
which suffices to obtained the first statement of the theorem. For the 
second statement, the symmetric matrix $ (Q_2 \otimes Q_1 ) \Gamma (Q_2 \otimes Q_1 )$ is estimated consistently by $(\^Q_2\otimes \^Q_1) \^\Gamma (\^Q_2\otimes \^Q_1)$ and so are its eigenvalues.
\qed

\subsection*{Proof of Proposition \ref{bura}}
We can notice that $\sqrt n \^Q_1 \^M \^Q_2 $ has the same asymptotic law than $\sqrt n Q_1 (\^M-M) Q_2$ whose asymptotic variance is consistently estimated by $[(\^Q_2\otimes \^Q_1) \^\Gamma (\^Q_2\otimes \^Q_1)]^+$ (see the proof of Proposition \ref{th1}). \qed

\subsection*{Proof of Proposition \ref{bootsir}}
Recall that $\^K_i= (X_i - \overline{X})(\Psi_i - \overline{\Psi})$, $K_i^*= w_i(\^K_i-\overline{\^K})$ and define $K_i = (X_i - \E[X])(\Psi_i - \E[\Psi])$. First note that, by Slutsky's theorem, $\sqrt n \ \overline{K^*}$ has the same asymptotic law than $n^{-1/2} \sum_{i=1}^n w_i(\^K_i-\E[K])$. Then we can develop   
\begin{align*}
&n^{-1/2} \sum_{i=1}^n w_i(\^K_i-\E[K]) \\
&=  n^{-1/2} \sum_{i=1}^n w_i((X_i-\E[X])(\Psi_i-\overline{\Psi})^T-\E[K] )+ (\E[X]-\overline{X}) n^{-1/2} \sum_{i=1}^n w_i(\Psi_i-\overline{\Psi})^T\\
&= n^{-1/2} \sum_{i=1}^n w_i(K_i-\E[K])+n^{-1/2} \sum_{i=1}^n w_i(X_i-\E[X])(\E[\Psi]-\overline{\Psi}) ^T \\
&\qquad \qquad \qquad \qquad \qquad \qquad \qquad \qquad \qquad +(\E[X]-\overline{X}) n^{-1/2} \sum_{i=1}^n w_i(\Psi_i-\overline{\Psi})^T.
 \end{align*} 
Checking a Lindeberg condition as bellow to ensure the weak convergence of $n^{-1/2} \sum_{i=1}^n w_i(X_i-\E[X])$ and $n^{-1/2} \sum_{i=1}^n w_i(\Psi_i-\overline{\Psi})^T$, and using the Slutsky's theorem we get conditionally a.s.  
 \begin{align*}
n^{1/2} \ \overline{K^*} = n^{-1/2} \sum_{i=1}^n w_i(K_i-\E[K])+O_{\P}(n^{-1/2}).
 \end{align*} 
We can apply the multidimensional version of the Lindeberg's central limit theorem (see for instance \cite{rao1976}, Corollary 18.2), provided that 
\begin{align*}
\frac 1 n \sum_{i=1}^n \E[\|\^V^{-1/2} w_i \xi_i \|^2 \mathds 1 _{\{\|\^V^{-1/2} w_i \xi_i \|>\nu n^{-1/2}\}} |\^P]\overset{\as}{\lr} 0,
\end{align*}
where $\xi_i = \vecv(K_i-\E[K])$ and $\^V= \frac 1 n \sum_{i=1}^n (\xi_i-\overline {\xi})(\xi_i-\overline {\xi})^T$. The above convergence is a consequence of the Lebesgue domination theorem which ensure that each term of the sum goes to $0$, afterwards we can conclude by the Cesaro's Lemma. Thus we have proved that conditionally a.s. 
\begin{align*}
n^{-1/2}\^V^{-1/2} \sum_{i=1}^n w_i\xi_i \overset{\dd}{\lr} \mathcal N (0,I),
\end{align*} 
and it remains to note that $\^V\overset{\as}{\r} V$ the variance of the limit in law of $\sqrt n(\^C - C)$ provided that $K$ has a finite order $2$ moment. For the second convergence, we note that conditionally a.s.
\begin{align*}
V^*-\^V &= \frac 1 n \sum_{i=1}^n (w_i^2-1) \xi_i \xi_i^T + o_{\P}(1),
\end{align*}
then by noting $v_i$ a coordinate of $\xi_i\xi_i^T$ we calculate
\begin{align*}
\E\left[\left(n^{-1}\sum_{i=1}^n (w_i^2-1)v_i\right)^2\right]=n^{-2}\E[(w_i^2-1)^2] \sum_{i=1}^n v_i^2
\end{align*}
which goes to $0$ a.s. provided that $K$ has a finite order $4$ moment. We conclude by using the Markov inequality to get that $V^*\overset{\P}{\r}\^V$ conditionally a.s..
\qed

\bibliographystyle{plain}
\bibliography{dimensionestimation}

\end{document}